# Lots of Aperiodic Sets of Tiles

Chaim Goodman-Strauss
Univ. Arkansas
strauss@uark.edu

**Abstract**

Aperiodic tiling — a form of complex global geometric structure arising through locally checkable, constant-time matching rules — has long been closely tied to a wide range of physical, information-theoretic, and foundational applications, but its study and use has been hindered by a lack of easily generated examples. Through readily generalized, robust techniques for controlling hierarchical structure, we increase the catalogue of explicit constructions of aperiodic sets of tiles hundreds-fold, in lots, easily assembled and configured from atomic subsets of 211 tiles, enforcing 25,380 distinct “domino” substitution tiling systems.

An **aperiodic set of tiles** is one that may be used to tile the plane, but only non-periodically — a form of complex global geometric structure arising through locally checkable, constant-time matching rules.

The very existence of aperiodic sets of tiles is implied by the undecidability of the “domino” (or “tiling”) problem, that no algorithm can ever decide whether any given set of tiles can be used to form a tiling of the plane. Hao Wang opened this discussion in 1961 [19], in the context of his work on one of the then-remaining open cases of Hilbert’s *Entscheidungsproblem* (“Is a given first order logical formula satisfiable?”). Wang conjectured that the domino problem is decidable, citing the self-evident implausibility of any existence of aperiodic sets of tiles! Fortuitously, Wang was incorrect and Robert Berger soon showed the domino problem undecidable [2], producing the first aperiodic set as a tool in his proof. Since this initial construction, about forty more aperiodic sets of tiles have been explicitly described, most found by mysterious art. These aperiodic sets of tiles have long been closely tied to a wide range of physical, information-theoretic, dynamical and foundational applications, in a range of geometric and combinatorial settings — see [6] for background, supporting and bibliographic material. However the study and use of aperiodicity through local rules has been hindered by a lack of easily generated examples.

Through readily generalized, robust techniques for controlling hierarchical structure, we increase the catalogue of explicit constructions of aperiodic sets of tiles hundreds-fold, in lots, easily assembled and configured from smaller atomic subsets, industrializing their production and flexibly enforcing a range of hierarchical, substitution tiling needs at reasonable cost, an example of control one might routinely expect. Such constructions may serve as scaffold for further applications and stimulate further development of the theory of matching rule tiling spaces.

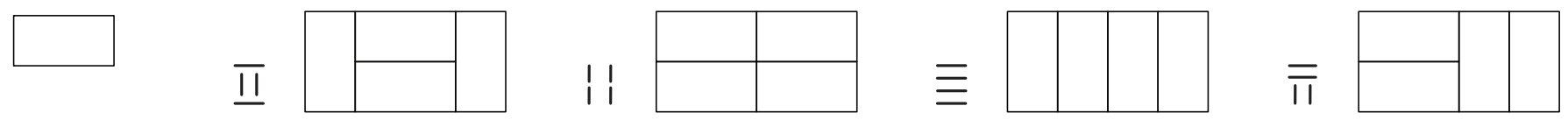

Figure 1: Four configurations of domino tiles, and suggestive notation for referring to them.

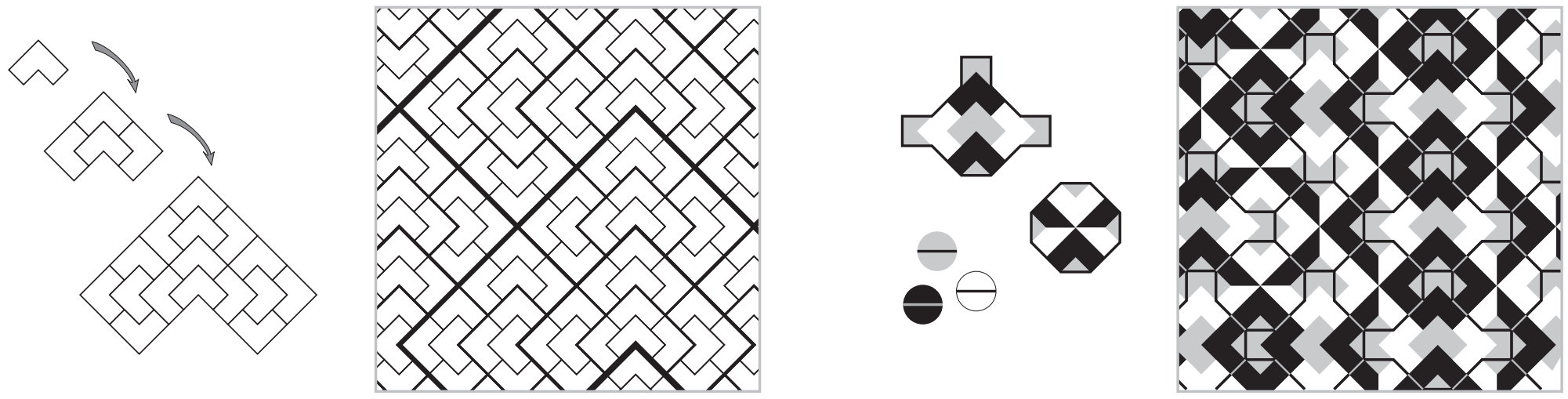

Figure 2: Enforcing substitution tilings with an aperiodic set of tiles: (Left) Non-periodic hierarchical tiling by "$L$-supertiles" defined by a tiling substitution rule (Section 1); however $L$-tiles admit many other tilings as well. At right, two tiles, the trilobite and crab [5], that *enforce* these hierarchical tilings: Each tiling by these two tiles *must* be a marked up hierarchical $L$-tiling. These two tiles *do* admit tilings (namely the marked up $L$-tilings) but only *non-periodic* ones (namely the marked up $L$-tilings and nothing else). Consequently the pair is an *aperiodic set of tiles*. We industrialize the production of such objects.

**Enforcing domino substitution tilings:** We follow the long thread from [2] onwards, constructing aperiodic sets of tiles that only admit hierarchical, and hence non-periodic, tilings (Figure 2), in our case on **domino** tiles, $2 \times 1$ rectangles

Our aperiodic sets "enforce" "substitution rules", as in Figure 2; we define these terms precisely in Section 1. Even the simplest non-trivial substitution rules on the domino are unexpectedly rich: of the four configurations, in Figure 1 (together with the symbols we'll use to name them) only three of them are enough to specify well-defined tiling substitution rules and hence tilings (see Figure 3). Each of the first three rules can be iterated in only one way. (The dynamics of the $\underline{\overline{\text{II}}}$ (or "table") tiling substitution system, at left of Figure 1, and of the $L$- (or "chair") substitution of Figure 2 were studied in [14].)

However the fourth rule is not yet well-defined: The domino tile has more symmetry than the fourth configuration, and so there is an ambiguity when we try to iterate, and this is more so when we try a second time. In order to give a well-defined rule we must give the specific motions that we are allowed to use to place each child supertile into its parent, so that we know which end is which as we iterate the rule. We address this by framing each supertile, with the markings of Section 4.1.

By specifying which "pieces" — specified by "atomic symbols" — of substitution rules we will allow (Figure 3), in other words, in which orientations we will allow children to be placed relative to their parents, we obtain a large number of distinct substitution tiling systems and their corresponding tiling spaces as (well-) defined from [9] onwards.

To each of these atomic symbols , we assign an "atomic set of tiles". We show that any union of atomic sets of tiles enforces the substitution rule that is described by the corresponding atomic symbols.

**Organization of the paper:** Much supporting, motivational, foundational and bibliographic material has been relegated to [6]. We define our terms in Section 1 and our markings in Section 2.

In Section 3 we construct a set of 4+23 tiles, $\mathcal{T}_1$ (Figure 8). Our atomic sets of tiles in $\mathcal{T}_1$ correspond to the three rules at left of Figure 1, individually, and all of the $\overline{\overline{\Pi}}$-rules, together, at right. As expressed in Theorem 1, there are nine aperiodic subsets of $\mathcal{T}_1$, unions of these atomic sets, enforcing the domino tiling substitution rules of Figure 1, together or individually, taking all the $\overline{\overline{\Pi}}$ tiling substitution rules as a group. The latter half of this section is somewhat *ad hoc* and may be easily passed over, though it does provide an example of a general technique used to work out atomic sets of tiles.

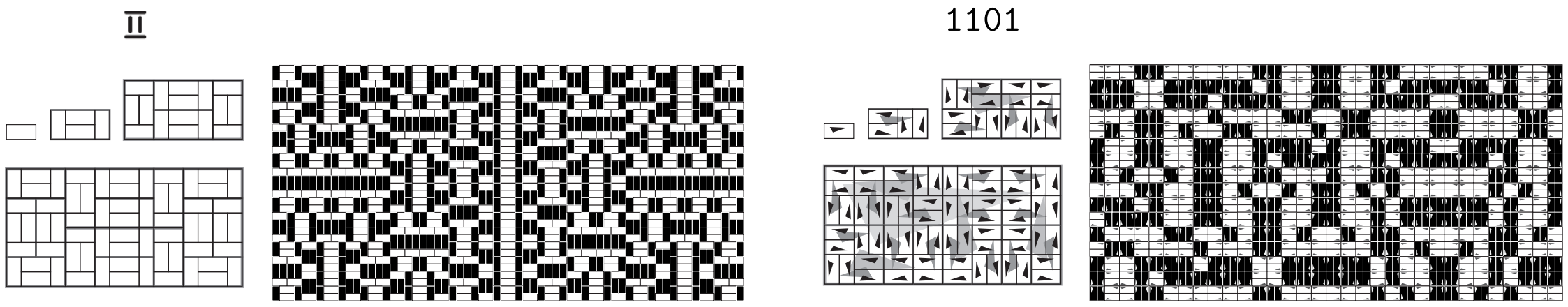

Figure 3: (Left) The "table" $\underline{\text{II}}$-tiling substitution [14], and a level 5 $\underline{\text{II}}$ supertile, with tiles marked by orientation. (Right) In order to apply a well-defined "domino" $\overline{\text{II}}$ tiling substitution rule, we must specify the isometries that place each tile into the supertile and child supertiles into parent ones — which way round is the placement of the tile? In the notation introduced in Section 4.2, the $\overline{\text{II}}$ tiling substitution rule shown at right above has code 1101.

In Section 4 we use a particular subset, denoted $\mathcal{T}_{\overline{\text{II}}}$, as the basis for a more refined set of 16+195 tiles, $\mathcal{T}_2$. Our atomic sets of tiles in $\mathcal{T}_2$ correspond to the manner in which each child domino may be placed into a parent $\overline{\text{II}}$ rule as suggested in Figure 17. As as made precise in Theorem 2, we may combine atomic subsets of $\mathcal{T}_2$ to enforce the $\overline{\text{II}}$ tiling substitution rules, together or individually. These atomic subsets exactly correspond to pieces of the structure of the $\overline{\text{II}}$ rules themselves.

In Section 5 we carefully count out 25,380 distinct matching rule tiling systems taking into account symmetries of the corresponding $\overline{\text{II}}$ substitution systems.[1] Of these, 128 are deterministic, with minimal substitution tiling spaces; the rest are composed from these minimal systems.

In order to verify, or not, that the corresponding tiling spaces are topologically distinct, we hope for the development of industrial strength computation of tiling space invariants, and a suitable theory of non-deterministically hierarchical matching rule tiling spaces.

In Section 6 we describe pairs of tiling substitutions that act upon a single tiling space. The corresponding tiling spaces are non-periodic yet have non-unique decomposition — each of these is a pair of distinct substitution rules that both act on the same tiling space (the action being in the sense of [9]). In this sense, these examples are not counter to the result of [18].

**Acknowlegements:** The sets of tiles $\mathcal{T}_{\underline{\text{II}}}$ and $\mathcal{T}_1$ were jointly discovered with Thomas Fernique, at his suggestion to enforce the table, $\underline{\text{II}}$ tiling substitution rule, while hosting me at Paris XIII. A few loose ends remained, and then several more, expanding into this current project, leaving many further questions to be considered. I am endebted to Thomas' collaboration and generosity, and to the hospitality of Le Laboratoire d'Informatique de Paris-Nord. Edmund Harriss and the reviewers have provided many helpful comments.

[1] Consider the ratio of the "size of a set of tiles" (specifying our geometric setting, definition of local rules, measure of complexity, etc.) to the "number of distinct tiling spaces we may enforce" with subsets. Berger's initial construction, famously, had a ratio of more than 20,000 (tiles to 1 tiling space). Robinson's tiles have a ratio of 6. Our set $\mathcal{T}_1$ has a respectable ratio of 3 (27 tiles to 9 tiling spaces), and there are many examples with a ratio of 2.The well-known "Einstein Problem" — *Is there an aperiodic monotile?* — asks for a ratio of 1 (and any solution tightly depends on specifying our setting and definition of complexity). For what it's worth, our second master set of tiles, $\mathcal{T}_2$ presents a substantial improvement in this metric (and in [6] we give an example of 84 tiles enforcing half a billion substitution systems!). It may be easy to settle: Must this ratio have infimum 0?

# 1 Definitions

We discuss foundational material more fully in [6]. Here let $\mathcal{I}$ be the set of Euclidean isometries, acting on $\mathbb{E}^2$ and let $\mathcal{A}^+$ be the set of distance expanding affine transformations of $\mathbb{E}^2$. Though we will illustrate and describe our tiles using colored markings, and the notion of local matching rules generalizes considerably, formally a **tile** $t$ is simply an unadorned, closed disk in $\mathbb{E}^2$. Our matching rules are enforced by restricting our definition of **configuration** to those sets $C$ of tiles with disjoint interiors.

The **support** $[\![C]\!]$ of a configuration $C$ is the union of the tiles, as point sets, within it, that is, $[\![C]\!] := \cup_{t \in C} t$. A **tiling** is a configuration with support $\mathbb{E}^2$.

In $\mathbb{E}^2$, a tiling $C$ is **periodic** if there exists some translation $v$ such that $C = C + v$; that is, for each tile $t \in C$, the tile $t + v$ is also in $C$.

A configuration $C$ is **admitted** by a set $\mathcal{T}$ of tiles if and only if each tile in $C$ is congruent to some tile in $\mathcal{T}$ (that is, $C = \{g_i t_i\}$ where each $g_i \in \mathcal{I}$ and each $t_i \in \mathcal{T}$). Let $\mathcal{C}(\mathcal{T})$ be all configurations admitted by $\mathcal{T}$ and let $\Sigma(\mathcal{T}) \subset \mathcal{C}(\mathcal{T})$ be all tilings admitted by $\mathcal{T}$.

A set of tiles is **aperiodic** if and only if it does admit a tiling, but does not admit any periodic tiling. In the Euclidean plane, this is satisfactory [8], but in general we must distinguish between *weakly aperiodic* sets of tiles, that admit only tilings with no compact fundamental domain (no co-compact symmetry) and somehow seem "for free" and more subtle *strongly aperiodic* sets of tiles, which destroy all infinite cyclic action (or all action).

We narrowly define a **tiling substitution rule** on a tile $t$ in a set $\mathcal{T}$ of tiles with "inflation factor" $s \in \mathcal{A}^+$, to be a configuration $\sigma(t) \in \mathcal{C}(\mathcal{T})$ such that $[\![\sigma(t)]\!] = st$; given a tiling substitution rule $\sigma$ on $t$, for any $g \in \mathcal{I}$, we define $\sigma(gt) := sgs^{-1}\sigma(t)$.

A finite set $\mathcal{S} = \{\sigma_i\}$ of tiling substitution rules $\sigma_i$ on tiles in some set $\mathcal{T}$ of tiles, at least one for each tile in $\mathcal{T}$, defines a **tiling substitution system** $(\mathcal{T}, \mathcal{S})$ (which is "deterministic" if there is exactly one rule per tile). We inductively define **supertiles** produced by $\mathcal{S}$:

- Each $\{t\}$, $t \in \mathcal{T}$, is a 0-level supertile.
- If $C$ is an $n$-level supertile, then each $gC$, $g \in \mathcal{I}$, is an $n$-level supertile.
- If $C$ is an $n$-level supertile, then $\displaystyle\bigcup_{gt \in C} \sigma(gt)$, each $\sigma \in \mathcal{S}$, is an $(n+1)$-level supertile.

Generally, for combinatorial substitution rules, we must take as axiomatic that somehow, each higher-level supertile is a well-defined configuration. Indeed, it remains an open question whether or not it is decidable if a given finite combinatorial substitution system is geometrically realizable! But here, because we assume that each $[\![\sigma(t)]\!] = st$ for a common $s \in \mathcal{A}^+$, by induction each supertile is in fact a well-defined configuration.

Note that as a configuration admitted by $\mathcal{T}$, $\sigma(t)$ has tiles of the form $g_i t_i, g_i \in \mathcal{I}, t_i \in \mathcal{T}$. Though often overlooked, in order to iterate substitution rules, we must specify these $g_i$, at least up to the symmetry of the supertiles (of all levels). In this paper, we will be concerned with the substitution rules defined by the configurations in Figure 1; the first three configurations have the same symmetry as the domino, and in order to define the $\overline{\underline{\mathrm{II}}}$, ¦¦ and ☰ substitution rules, the specific isometries used to place a given domino need only be specified up to such a symmetry. However, every higher level $\overline{\pi}$ supertile has trivial symmetry, and in order to be able to well-definedly iterate a $\overline{\pi}$ substitution rule, we must give the actual isometries used to place each domino.

We take various sets of these rules to define non-deterministic substitution rules, such as {❙❙, ☰}, or various collections of the deterministic ╥╥ rules.

The **substitution tilings** defined by a given tiling substitution system $(\mathcal{T}, \mathcal{S})$ are those tilings $C \in \Sigma(\mathcal{T})$ such that for each finite subset $C_1 \subset C$, there exists some supertile $C_2$ and an isometry $g$ such so that $gC_1 \subset C_2$. The **hierarchical tilings** defined by $(\mathcal{T}, \mathcal{S})$ are the tilings $C \in \Sigma(\mathcal{T})$ such that for each $t \in C$, there exists a sequence $t = C_0 \subset \ldots \subset C_n \subset \ldots \subset C$, each $C_n$ an $n$-th level supertile. It is not difficult to show that (a) every substitution tiling is a hierarchical tiling and (b) under any translation invariant probability measure, almost every hierarchical tiling is a substitution tiling.

(Almost every such hierarchy of supertiles covers the entire plane, though still uncountably many will cover only a portion of the plane. A substitution tiling system "forces the border" [17] if and only if the substitution tilings made of up of more than one infinite hierarchy (each partially covering the plane) are each determined by any one of these hierarchies. )

Otherwise, as in [15] onwards, in a hierarchical tiling, there may be an "infinite fault line" between these hierarchies, with no coordination across the fault. However, in any tiling space, as defined in [13], the subset of tilings with infinite faultlines can only be of measure zero, and so is often disregarded.)

(Non-deterministic substitution tilings, such as those we consider here, do give well-defined tiling spaces, in precisely the sense of [13] onwards, but these are not minimal. A foundation for the study of such spaces appears in [16], and much further development may be hoped for.)

**Matching rules enforcing hierarchical tiling:** Many different definitions of "enforcing" a tiling substitution system appear in the literature[2], but all of them imply an almost-everywhere well-defined, onto, locally-checkable map from a matching rule tiling space to a substitution tiling space, or equivalently, up to measure one, a hierarchical tiling space.

Our definition here is aligned with the specific structure of our inductive proof of enforcement, following [2] onwards and will satisfy whatever reasonable demands we make of it.

Given sets of tiles $\mathcal{T}', \mathcal{T}$ a **local map**, or **local derivation** $\Phi : \Sigma(\mathcal{T}') \to \Sigma(\mathcal{T})$ is defined by specifying a finite set of distinct finite configurations $\{C_i\} \subset \mathcal{C}(\mathcal{T}')$ and a corresponding set of finite configurations $\{\Phi(C_i)\} \subset \mathcal{C}(\mathcal{T})$, satisfying:

- Each tiling $\tau \in \Sigma(\mathcal{T}')$ is the union of configurations of the form $gC$, $g \in \mathcal{I}$, $C \in \{C_i\}$.
- $\Phi(\tau) := \bigcup_{gC \subset T} g\Phi(C)$ is a tiling in $\Sigma(\mathcal{T})$.

We say that a set $\mathcal{T}'$ of tiles **enforces** a substitution system $(\mathcal{T}, \mathcal{S})$ if and only if there is a well-defined local map $\Phi : \Sigma(\mathcal{T}') \to \Sigma(\mathcal{T})$ such that each $\Phi(\tau)$ is a hierarchical tiling under $(\mathcal{T}, \mathcal{S})$.

(In fact, we can strengthen our definition of enforcement, and require that each $\Phi(\tau)$ is a substitution tiling – that is we can avoid "slipping along infinite fault lines" — by adding a sidedness to our marking +, − below, as our other markings have, with an increase in the size of our tilesets. We touch on this and other variations in [6].)

[2]Other well-established, equivalent-up-to-measure-one, definitions of "enforcing tiling substitution system" include: requiring $\Phi$ to be a bijection from all but a measure zero subset of $\Sigma(\mathcal{T}')$ (measured in any translation invariant Borel probabability measure) to the substitution tilings of $(\mathcal{T}, \mathcal{S})$ (cf. [12]). Or we may define preimages of supertiles, and require that almost every tile in every tiling in $\Sigma(\mathcal{T}')$ lies in a unique hierarchy of such preimages (cf. [4]). Or we may require that the tiles in $\mathcal{T}'$ are decorated versions of the tiles in $\mathcal{T}$ and the map $\Phi$ that simply removes the decorations is a bijection from $\Sigma(\mathcal{T}')$ to the set of substitution tilings. We do not discuss the meaning of these variations here, and our construction is transparent to them.

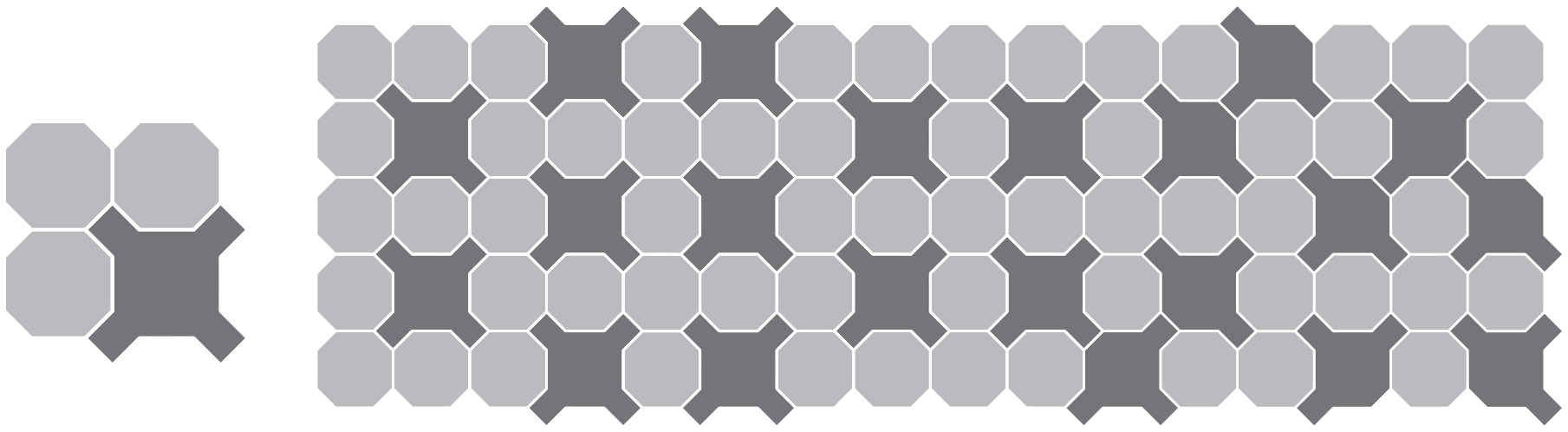

Figure 4: At left, cornered and uncornered tiles, following Robinson [15]. The tiles in any tiling by these lie on a square lattice; each vertex of the lattice meets one cornered tile and three uncornered tiles, and each cornered tile must be in at least a row or a column of cornered tiles. As at the right end of the illustration, the "trilobite" tile has richer behavior and, as in [5], may allow a more efficient implementation of our construction here and more flexibly admit more complex structures.

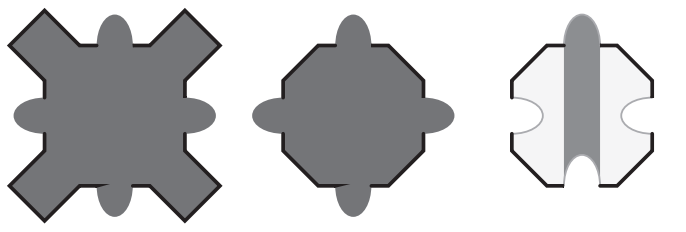

Figure 5: Outward and crossing tiles.

## 2 The underlying markings and tiles

We *define* our matching rules, as above, geometrically: our tiles are unadorned closed topological disks and we require that these have disjoint interiors in any well-formed configuration. We will *describe* and *enumerate* our matching rules and tilings, encoded as below, as various strings on edges, and we will *illustrate* them with colors, as in Figure 6. We may adopt other formalisms if we wish instead, as for example Wang tiles in [15].

Our basic construction is, as will be easily recognized, derived from the Robinson tiles [15], and as there, our tiles are modified squares, of two basic kinds: "cornered" and "uncornered" which may only fit together as in Figure 4.

Our markings are all "directed", pointing in or out, and some are "sided", balanced right or left, captured in the notation below. A tile is defined, up to congruence, by specifying whether it is cornered or uncornered, and its four edge markings, up to cyclic ordering, or reversing the ordering and reflecting the markings. We always presume our tiles are denoted in the (*ad hoc* as useful) normal forms we develop for referring to them.

As indicated in Figure 5, we will be working with "outward" tiles (which may be cornered or uncornered), with all four markings "directed outward", and uncornered "crossing" tiles, with three markings "inward", and one outward marking matching the opposite inward one, as described below

The markings themselves are illustrated in Figure 6 and the basic structure they enforce in Figure 7.

**Encoding the markings symbolically:** Our edge markings are of the form

$$\mathtt{m} = (\mathtt{abcd})_{\mathtt{m}} \in \{\mathtt{+},\mathtt{-}\} \times \{\mathtt{0},\mathtt{+},\mathtt{-}\} \times \{0,1,2,3\} \times \{0,1,2,3\}$$

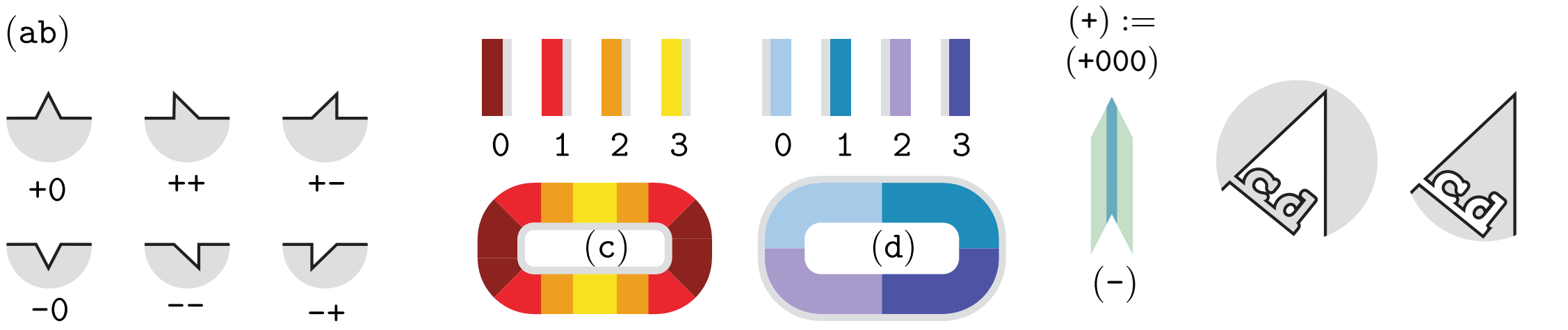

Figure 6: Our basic edge markings, (abcd) have four channels, illustrated here. The (ab) channel are precisely the Robinson tile markings [15]. Tiles may meet along an edge if and only if the (ab) channels sum to 00 and the (cd) channels are equivalent. Reflecting a tile takes (b) to (-b). The special markings ($\pm$000) are denoted simply (+), (-). $\mathbb{Z}_2 \oplus \mathbb{Z}_2$ acts upon the d-channel by nim-addition (§4.1), just as it acts upon the symmetries of the domino. Though we illustrate our constructions with colored markings, and define them as combinatorial notations, we can encode these as simple topological rules, our tiles being unmarked closed topological disks, with combinatorics encoded within the geometry, as at right. More generally, as discussed in [6], we may refine and adjust our markings, for different ends, with varying control of our tiling spaces.(Specifying the sidedness of (000), for example, we can coordinate infinite supertiles across "infinite fault lines", but with a correspondingly greater number of tiles.)

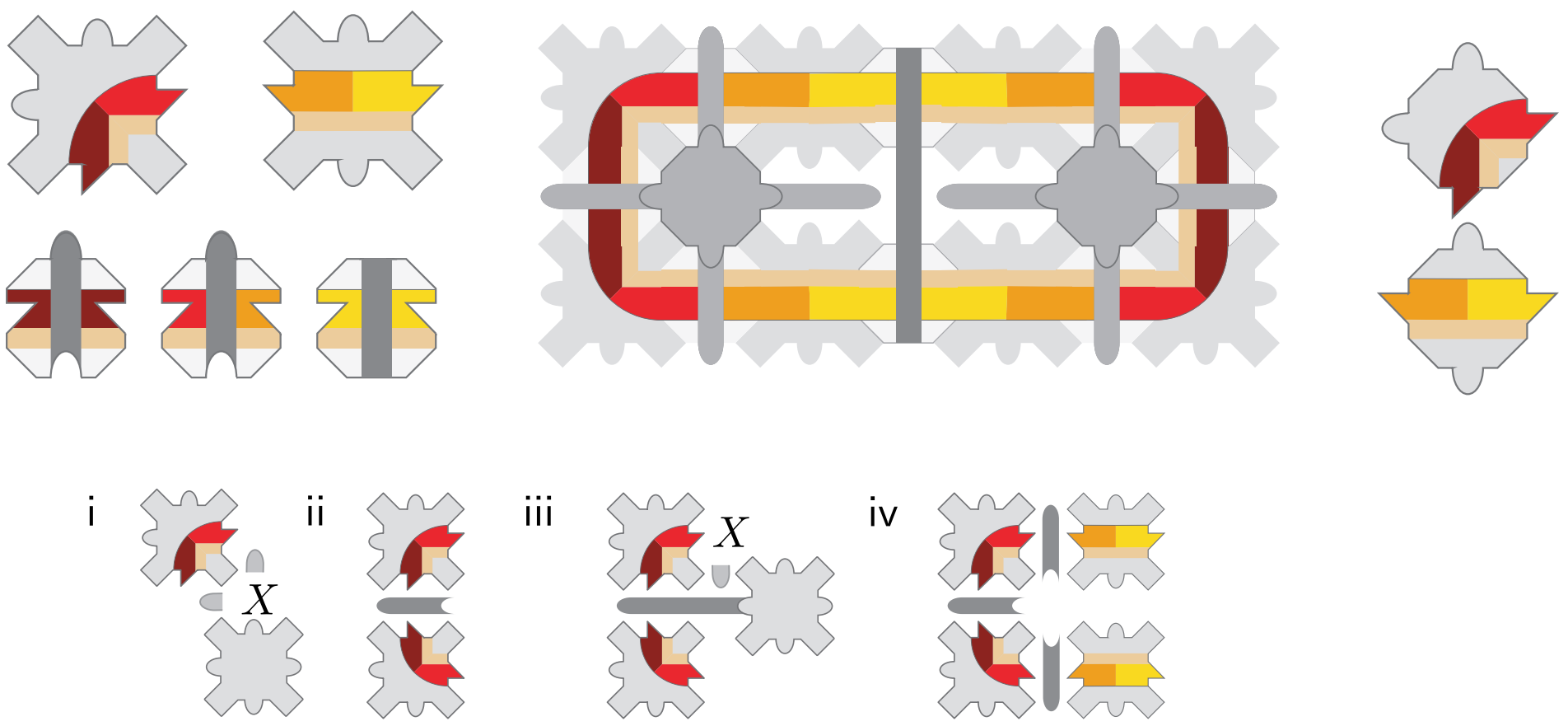

Figure 7: The tiles in $\mathcal{T}_1$ will be based on those at top left and right, with the remaining markings given in Figure 9. As shown at middle left, upon each crossing tile with horizontal (c) channel markings 0, 1 or 2, the vertical marking must be "up", i.e. (a) channel marking + on top of the tile and - on bottom. On those marked 3, the vertical markings may be up or down. (The active reader is encouraged to write this notation on the figure, to facilitate further calculations throughout.) Consequently, in any tiling by tiles with at least these markings, the cornered tiles must each lie within the $3 \times 7$ "domino block" shown at middle top, as we show below: (i) If cornered tiles are placed as at lower left, no tile may be placed at the location marked $X$. (ii) Consequently, copies of the tile at upper left must occur in pairs. (iii) Similarly, no configuration includes the tile shown and a tile at $X$. (iv) Thus, each cornered tile lies within a quadruple, well-defining some domino block. In Section 3.3, in the interior of each domino, we will place the uncornered outward tiles shown at right, controlling their possible placements as "atomic sets" of allowed "key tiles", where signed crossings meet. In Section 4 we gain finer control through the use of our d-channel markings.

($\mathtt{d}$ will be suppressed in $\mathcal{T}_1$).

Edge markings $\mathtt{m}$, $\mathtt{m}'$ may be fitted together (that is, "match") if and only if $(\mathtt{ab})_{\mathtt{m}} = \mathtt{-}(\mathtt{ab})_{\mathtt{m}'}$, and $(\mathtt{cd})_{\mathtt{m}} = (\mathtt{cd})_{\mathtt{m}'}$ (so for example, (+-31) and (-+31) match). A marking $\mathtt{m}'$ is a "reflection" of a marking $\mathtt{m}$ if and only if $(\mathtt{b})_{\mathtt{m}} = (-\mathtt{b})_{\mathtt{m}'}$, and $(\mathtt{a}, \mathtt{cd})_{\mathtt{m}} = (\mathtt{a}, \mathtt{cd})_{\mathtt{m}'}$ (so for example, (+-31) and (++31) are reflections of one another).

In the illustrations of markings, described in Figure 6, we encode these properties as geometry. With no ambiguity in the markings we will use, we may abbreviate (+000) and (-000) as (+) and (-), and if d is unneeded (as in the first set of tiles), we suppress it, writing (+-3) for (+-30).

We naturally refer to the (a, b, c or d) channel of a marking m, meaning $(\mathtt{a})_{\mathtt{m}}$, $(\mathtt{b})_{\mathtt{m}}$, etc. The (a) channel is the "direction" of a marking, the (b) channel is its sidedness. The (c) channel is used to encode the structure of $\Sigma(\mathcal{T}_1)$. The d-channel further encodes $\Sigma(\mathcal{T}_2)$, which factors to $\Sigma(\mathcal{T}_{\text{╥}})$ by suppressing d.

The $\mathbb{Z}_2 \oplus \mathbb{Z}_2$ symmetry of the domino naturally operates on symbols in the d-channel, by "nim-additon", in which bits are added independently mod 2 as described further in Section 4.1.

# 3 Nine aperiodic subsets of twenty seven tiles.

Our goal in this section is to define the terms of the following theorem, and provide its proof:

**Theorem 1** *The set of tiles*

| | | | |
|---|---|---|---|
| $\mathcal{T}_{\text{Ⅱ}}$ | *enforces the* | Ⅱ | *tiling substitution system;* |
| $\mathcal{T}_{\text{¦¦}}$ | *enforces the* | ¦¦ | *tiling substitution system;* |
| $\mathcal{T}_{\text{≣}}$ | *enforces the* | ≣ | *tiling substitution system;* |
| $\mathcal{T}_{\text{╥}}$ | *enforces the* | {¦¦, ≣, ╥} | *tiling substitution system;* |
| $\mathcal{T}_{\text{Ⅱ}} \cup \mathcal{T}_{\text{¦¦}}$ | *enforces the* | {Ⅱ, ¦¦} | *tiling substitution system;* |
| $\mathcal{T}_{\text{Ⅱ}} \cup \mathcal{T}_{\text{≣}}$ | *enforces the* | {Ⅱ, ≣} | *tiling substitution system;* |
| $\mathcal{T}_{\text{¦¦}} \cup \mathcal{T}_{\text{≣}}$ | *enforces the* | {¦¦, ≣, ╥} | *tiling substitution system;* |
| $\mathcal{T}_{\text{Ⅱ}} \cup \mathcal{T}_{\text{¦¦}} \cup \mathcal{T}_{\text{≣}}$ | *enforces the* | {Ⅱ, ¦¦, ≣} | *tiling substitution system; and* |
| $\mathcal{T}_1$ | *enforces the* | {Ⅱ, ¦¦, ≣, ╥} | *tiling substitution system.* |

*Each of these sets of tiles is aperiodic.*

*We further observe that for each of the above sets, its union with $\mathcal{U}_1$ enforces just exactly what the set itself does on its own. No proper subset of $\mathcal{T}_{\text{Ⅱ}}$, $\mathcal{T}_{\text{¦¦}}$ or $\mathcal{T}_{\text{≣}}$ admits a tiling.*

In Section 3.1 we outline the structure of the proof, which is essentially the same as that of every aperiodic hierarchical tiling from [2] onwards.

In Section 3.2 we formally define a set of 27 tiles $\mathcal{T}_1$, for which various subsets enforce domino substitution tiling systems as claimed in Theorem 1.

In Section 3.3, we give the underlying combinatorial structure. In essence, our tiles must lie in well-formed marked supertiles, and we can control which ones by what crossings we allow, in the form of "key tiles", where the signed markings of a marked supertile meets those of its parent and grandparent (see Figure 8 and Figure 15).

Only then do we work out what what we need to include and exclude from our tile sets $\mathcal{T}_{\text{Ⅱ}}$ (16 tiles), $\mathcal{T}_{\text{¦¦}}$ (17 tiles),$\mathcal{T}_{\text{≣}}$ (18) and $\mathcal{T}_{\text{╥}}$ (26). The details are tedious and *ad hoc*, but this serves as an example of iterating substitutions on marked supertiles, a useful and general technique for working out such sets. In Section 3.5 we further check how these sets are related to one another, completing the definition of the terms of Theorem 1 and the details of its proof. Once accepting that $\mathcal{T}_{\overline{\overline{\shortparallel}}}$, particularly, behaves as stated, one may safely skip forward to Section 4, and its Theorem 2 which has a more systematic proof.

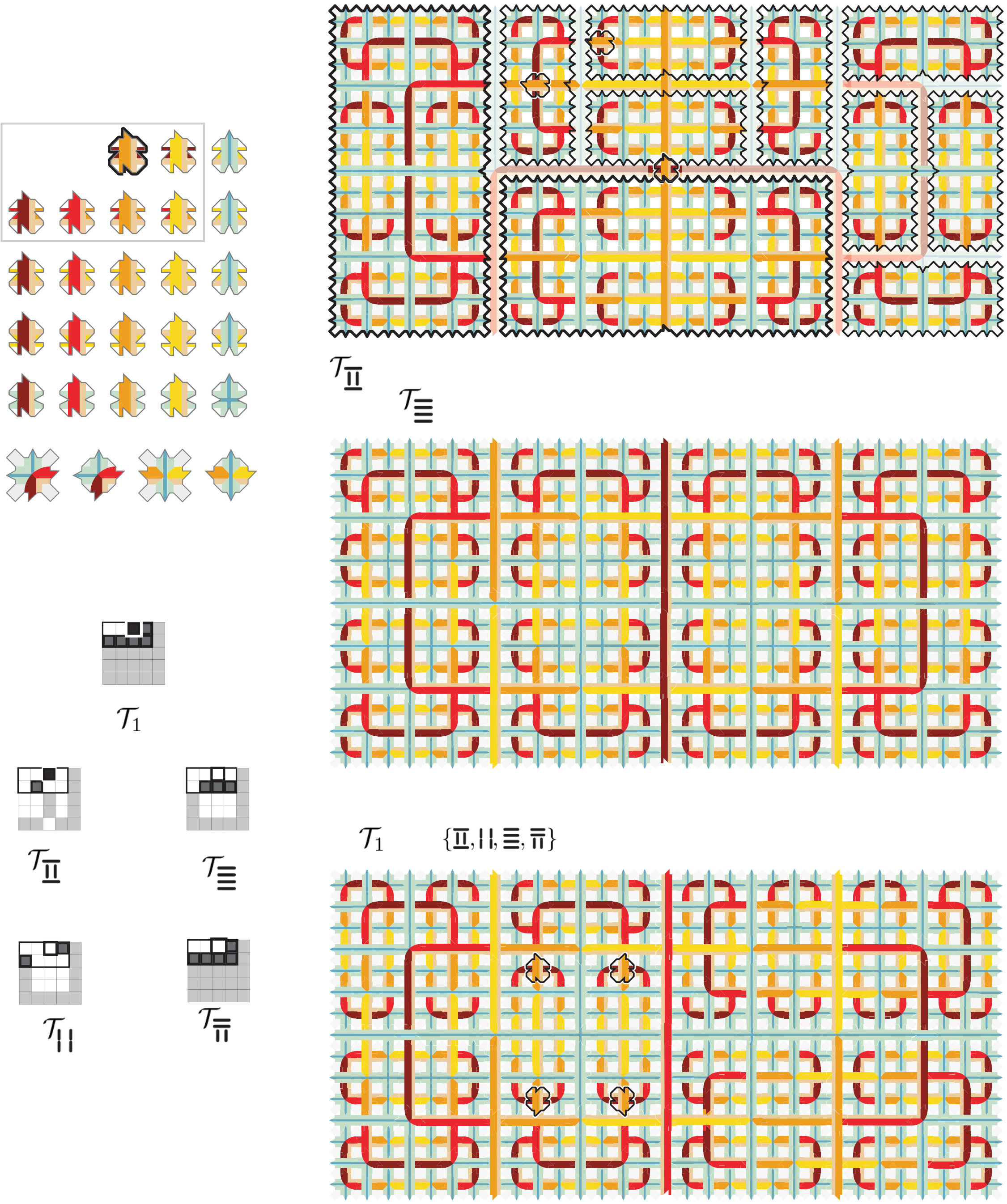

Figure 8: At left top we show the tiles of $\mathcal{T}_1$, boxing its “key tiles” where signed markings meet. Four of its aperiodic subsets, enforcing different domino tiling substitution systems are indicated below. We control which ╦ substitutions are allowed simply by which tiles we include or not. For example, the marked tile in the table at top left, is needed to form Ⅱ supertiles, as at top right and bottom right (in particular it is required for a `U` block to be within a parent `I` block). Without the marked tile we can enforce at most {¦¦, ≣, ╦}. (Top right) A marked Ⅱ supertile admitted by the 16 tiles of $\mathcal{T}_{Ⅱ}$. (Middle right) A marked ¦¦ supertile admitted by the 17 tiles in $\mathcal{T}_{¦¦}$. The marked tiles are the same, both lying at the center of an `H` within an `I`. (Bottom right) A {Ⅱ, ¦¦, ≣, ╦} supertile in the 27 tiles of $\mathcal{T}_1$.

### 3.1 Outline of the proof of Theorem 1

Our strategy is fairly standard, similar to each construction of matching rules enforcing hierarchical tilings from [2] onwards. Here and more so in the proof of Theorem 2, we industrialize our arguments, handling many standardized cases at once.

As usual in constructing aperiodic sets of tiles enforcing hierarchical structure, we show that each set of tiles enforces some specified tiling substitution system, each tile in each tiling lying in a well-defined, unique hierarchy of larger and larger "marked supertiles", configurations mapped under a natural local decomposition to supertiles in the tiling substitution system (by simply erasing the markings and adjusting the boundaries).[3] Consequently, the set of tiles actually admits a tiling at all (because because it admits configurations covering arbitrarily large disks — this standard argument appears as The Extension Theorem, Theorem 3.8.1 in [8] and is often cited as "by compactness" or "by Koenig's Lemma") and is aperiodic (as each tile lies in a unique marked supertile of each size).

We accomplish this by showing inductively, that every tiling by $k$-level marked supertiles is also a tiling by $(k+1)$-level marked supertiles.

As in Figure 8, we may restrict or allow which supertiles may be formed in the hierarchy by excluding or including particular tiles. In Section 3.4, in order to verify which tiles are actually needed, we let substitution act on the marked tiles and supertiles themselves, enumerating tiles as we go along.

The erasing mapping from the matching rule tiling space is well-defined onto the corresponding hierarchical tiling space (and almost everywhere well-defined onto substitution tiling spaces).

Each of the sets of tiles will enforce the tiling substitution systems specified in the theorem.

One note: The hierarchies of the ¦¦ tiling substitution are combinatorially those of the Robinson tiling, yet the two structures are not mutually locally derivable — the scalings are incompatible and the supertiles cannot be even quasi-isometric.

### 3.2 Defining $\mathcal{T}_1$

The first set, illustrated in Figure 9 denoted $\mathcal{T}_1$, consists of 27 tiles, all with `d` unused and suppressed. $\mathcal{T}_1$ has four outward tiles, which we denote $\mathcal{T}_{+}$, which may be uncornered or cornered and markings `(+)(+)(++1)(+-0)` and `(+)(+-2)(+)(++3)`, illustrated at top left of Figure 9.

Our twenty-three crossing tiles, $\mathcal{T}_{\mathtt{hv}}$ have edge markings of the forms shown in Figure 9. For simplicity, as in the table in at middle right of the figure, we denote these markings [`hv`] where $\mathtt{h} = \mathtt{-}, \mathtt{0}, \mathtt{1}, \mathtt{3}, \mathtt{-3}$ denotes the tile's row in the table (and is based off of the west markings); and $\mathtt{v} = \mathtt{+}, \mathtt{0}, \mathtt{1}, \mathtt{2}, \mathtt{3}$ denotes its column. Note that the tiles [`00`], [`01`] are not included in $\mathcal{T}_{\mathtt{hv}}$, and hence not in $\mathcal{T}_1 = \mathcal{T}_{+} \cup \mathcal{T}_{\mathtt{hv}}$.

The orientation of the vertical marking relative to the horizontal ones are exactly as shown

[3] As in [9], this extends naturally to a surjection from our matching rule tiling spaces onto the hierarchical ones, which contain substutition tilings as a subset of measure one. This map is one-to-one, except on a measure zero set, the hierarchical tilings with "infinite fault lines".

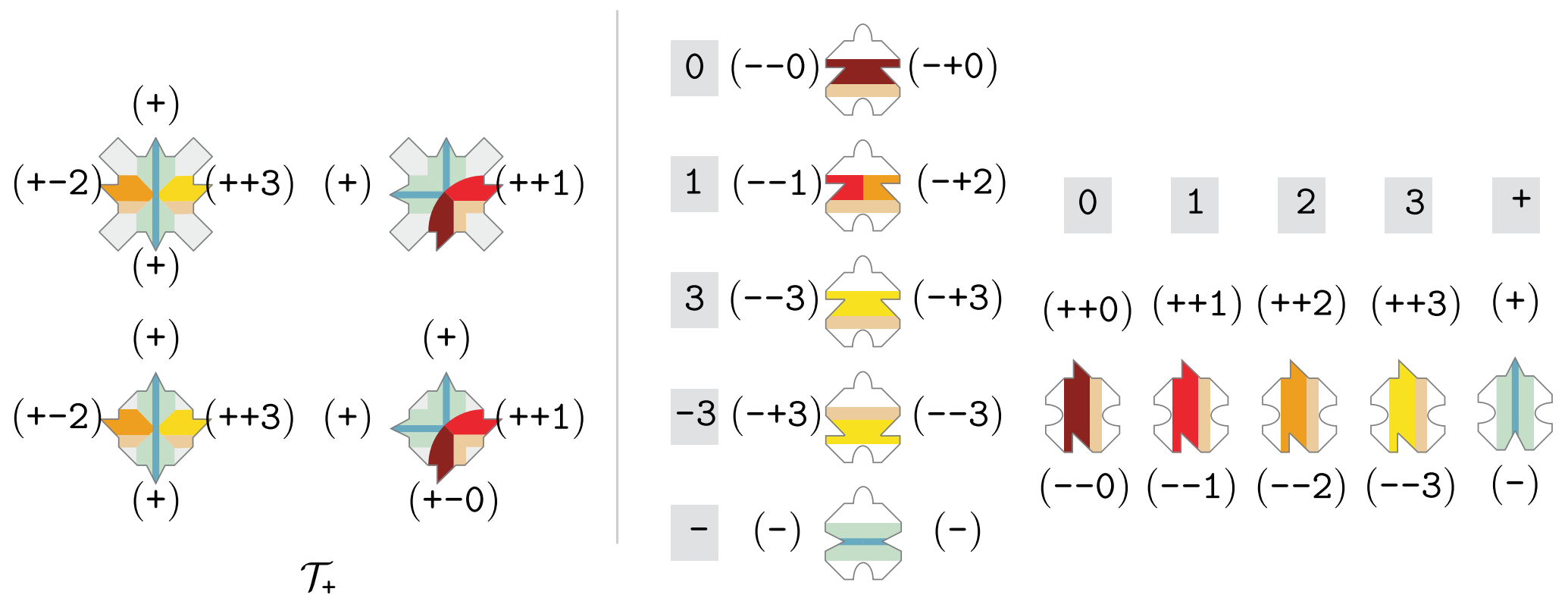

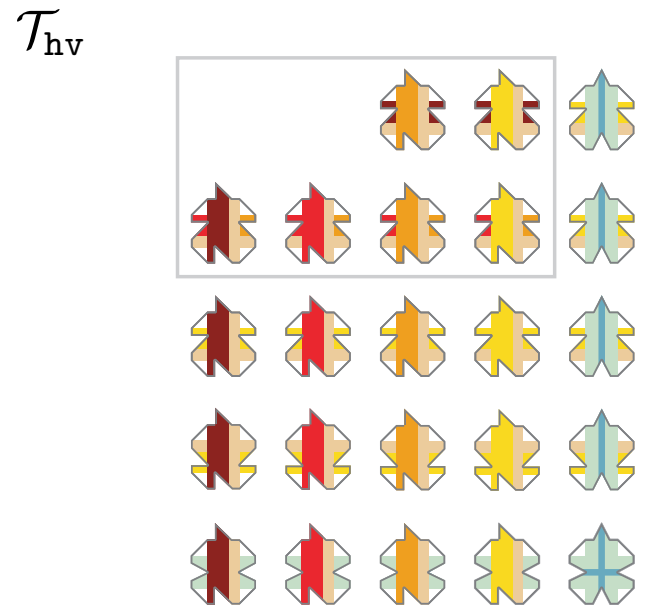

| h \ v | 0 | 1 | 2 | 3 | + |
|---|---|---|---|---|---|
| 0 | | | [02] | [03] | [0+] |
| 1 | [10] | [11] | [12] | [13] | [1+] |
| 3 | [30] | [31] | [32] | [33] | [3+] |
| -3 | [-30] | [-31] | [-32] | [-33] | [-3+] |
| - | [-0] | [-1] | [-2] | [-3] | [-+] |

Figure 9: We explicitly define the markings on our first set $\mathcal{T}_1$ of twenty-seven tiles, all with d suppressed. There are four outward tiles $\mathcal{T}_+$, shown at upper left, and twenty-three crossing tiles $\mathcal{T}_{\mathtt{hv}}$, shown at lower left (with six "key tiles" boxed). The markings making up the crossing tiles are illustrated in the upper right. At lower right, we give compressed notation to refer to the crossing tiles.

— which due to symmetry makes no difference except for the tiles with horizontal marking (--1)(-+2) and a vertical marking other than (+). In Figure 10 we see the necessity of this restriction.

We formally define:

We define

$\mathcal{T}_1 := \mathcal{T}_+ \cup \mathcal{T}_{\mathtt{hv}}$, where $\mathcal{T}_+$ is our set of outward tiles and

$$\mathcal{T}_{\mathtt{hv}} := [(0,1,3,\text{-}3,-)(0,1,2,3,+)] \setminus \{[00],[01]\}$$

$$= \{[02],[03],[10],[11],[12],[13],$$
$$[30],[\text{-}30],[31],[\text{-}31],[32],[\text{-}32],[33],[\text{-}33],$$
$$[\text{-}0],[\text{-}1],[\text{-}2],[\text{-}3],[\text{-}+],[0+],[1+],[3+],[\text{-}3+]\}$$

Define "key tiles", where signed markings cross, and the rest:

$$\mathcal{K}_1 := [(0,1)(0,1,2,3)] \setminus \{[00],[01]\} = [(0,1)(0,1,2,3)] \cap \mathcal{T}_{\mathtt{hv}}$$

$$\mathcal{U}_1 := \mathcal{T}_{\mathtt{hv}} \setminus \mathcal{K}_1$$

## 3.3 The underlying combinatorial structure

Consider any subset $\mathcal{T}$ of our twenty-seven tiles $\mathcal{T}_1$: If $\mathcal{T}$ is to admit any tiling at all, it must include one of the cornered outward tiles, and every tile in any such tiling is a cornered outward tile or shares a corner with one. As in Figure 7, these cornered tiles can only lie within well-defined $3 \times 7$ "domino blocks".

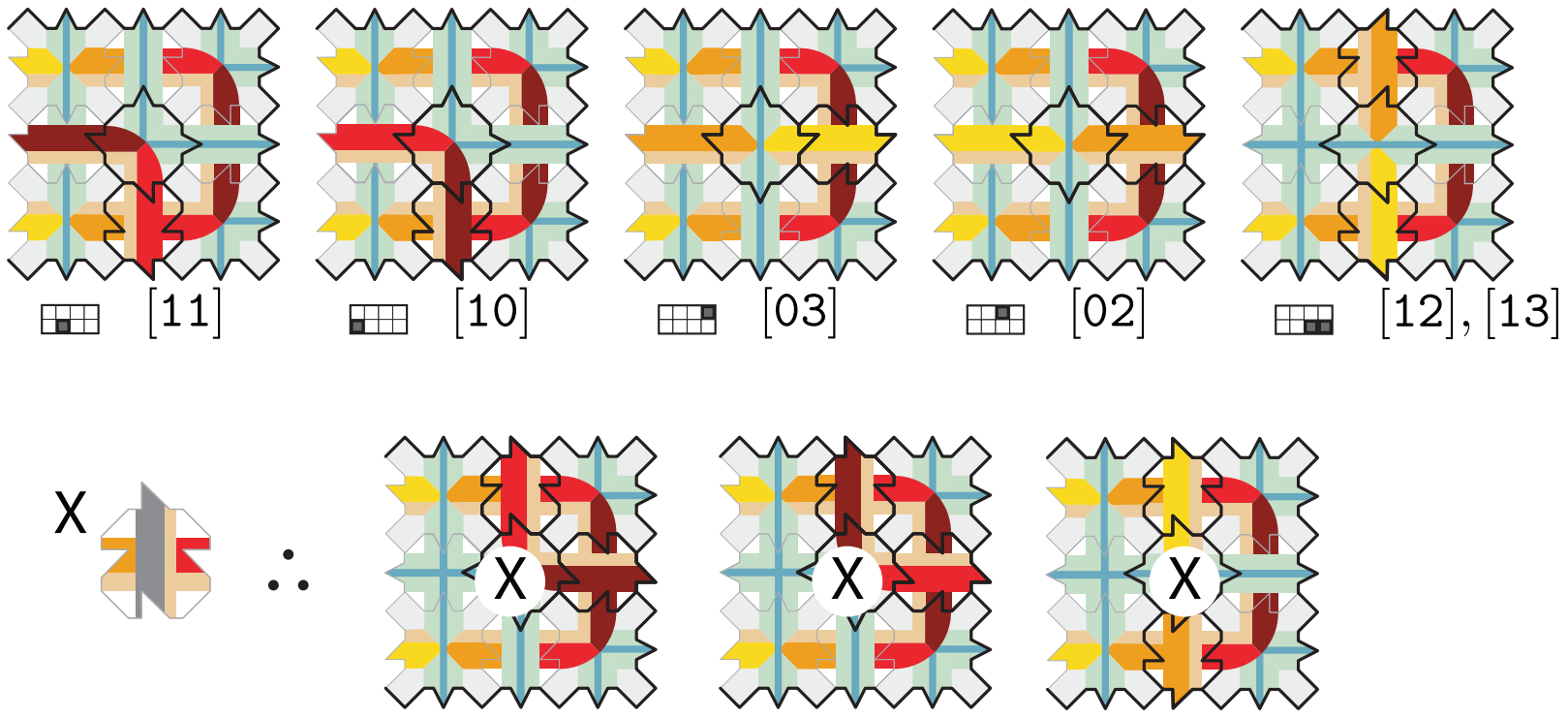

Figure 10: Cornered outward tiles can only lie in half domino blocks, specified by the type and orientation of the uncornered outward tile at its the center. There are only five possibilities, because of the asymmetry of the horizontal marking in the tile at bottom left. For each of the possible half blocks, one or two key tiles, outlined in the figure, are essential. These "atomic sets" are a signature for each type of block. (The rest of the crossing tiles, that are *not* highlighted, are more common and relegated to $\mathcal{U}_1$.)

In Figure 10, we consider the possible ends of these domino blocks, depending on the particular type and orientation of the uncornered outward tile at their centers. Because of the asymmetry of the $(1)-(2)$ horizontal markings in the $[1*]$ tiles, only five of these half dominos are possible,

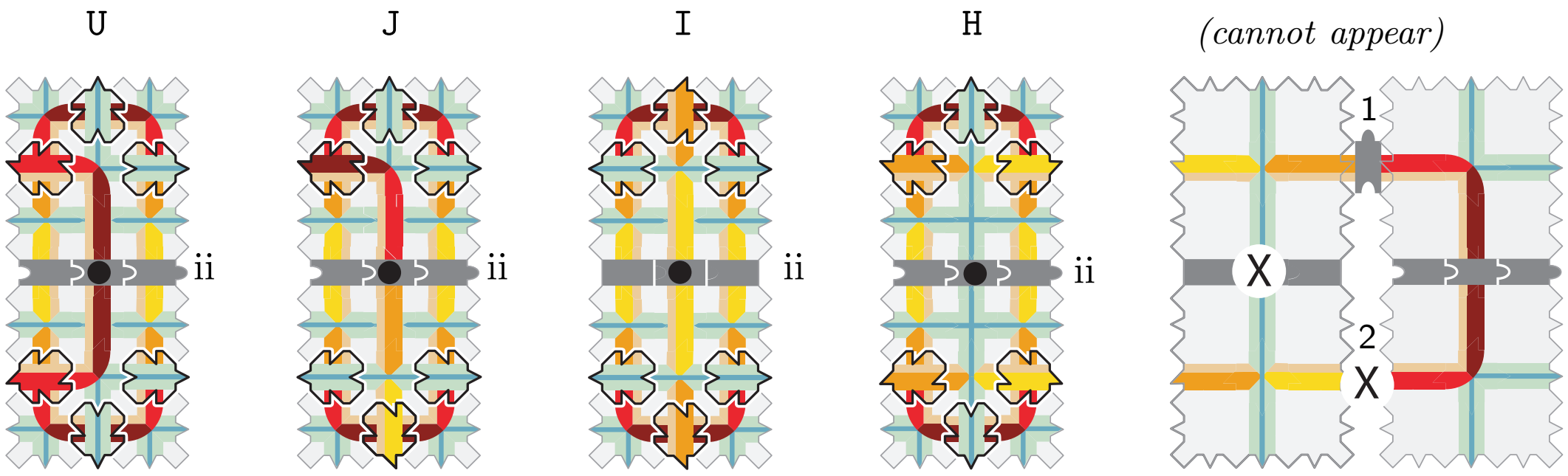

Figure 11: The tiles in $\mathcal{T}_1$ admit four $(2 \times 4)$- blocks, U, J, I, H — The markings, and most of the tiles, of each block are determined, leaving only a central band to be specified to fully define a block. Each block has matching markings on either end, their direction specified for U, J and H blocks as shown. . The sets of tiles $\mathcal{T}_{\mathtt{U}}$, $\mathcal{T}_{\mathtt{J}}$, $\mathcal{T}_{\mathtt{I}}$, and $\mathcal{T}_{\mathtt{H}}$ that they require are outlined above, and are given in compressed graphic form at left in Figure 12. A fifth block can be formed, but cannot appear in any tiling: As indicated at right, at (1), the marking (+−1) can only meet a marking (++0) of a block of type U. But then no tile may fit at (2).

and in Figure 11 we see that these can be united into exactly four possible $2 \times 4$ blocks, denoted H, I, J and U, up to just one remaining (c) marking, presumably from some higher level block, to be determined when useful to do so. In each case, there are a few tiles that are essential in order to assemble each type of block. We define:

| | | key tiles | the rest |
|---|---|---|---|
| 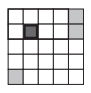 | $\mathcal{T}_{\mathtt{U}} :=$ | {[11], | [−0], [0+], [1+]} |
| 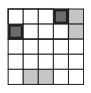 | $\mathcal{T}_{\mathtt{J}} :=$ | {[03], [10], | [−1], [−2], [0+], [1+]} |
| 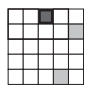 | $\mathcal{T}_{\mathtt{I}} :=$ | {[02], | [−3], [1+]} |
| 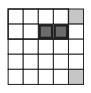 | $\mathcal{T}_{\mathtt{H}} :=$ | {[12], [13], | [0+], [−+]} |

For example, if a subset $\mathcal{T}$ admits a tiling that includes a copy of block J, then the tiles [01] and [30] must be in $\mathcal{T}$ — that is, if they are not, this block is in no tiling admitted by $\mathcal{T}$.

### 3.4 Deriving the sets $\mathcal{T}_{\text{II}}, \mathcal{T}_{\text{¦¦}}, \mathcal{T}_{\equiv}$ and $\mathcal{T}_{\overline{\text{II}}}$

The remaining details are fairly tedious and less illuminating than those of $\mathcal{T}_2$. One may safely skip the remainder of this section and proceed to Section 4, accepting that $\mathcal{T}_{\overline{\text{II}}}$ admits, and only admits U, H and J blocks.

We will calculate precisely which matching rule tiles are needed to allow, and only allow, particular blocks, by substituting on marked blocks, beginning with a generic $\mathtt{X}(x)$, and then continuing on any new blocks that arise, tallying up needed tiles as we proceed.

In Figure 12 we give indices for referring to the markings at various positions in any of our blocks, with respect to some arbitrary orientation. Note that the orientation (that is, the a and b parts) is determined in all the markings *i*, *iii*-*vii*. The direction, that is, the (a) channel, of the markings *viii* and *ix* is determined, and nothing about the marking *ii* is determined. A block will have a pair of matching markings *ii* (top) and -*ii* (bottom).

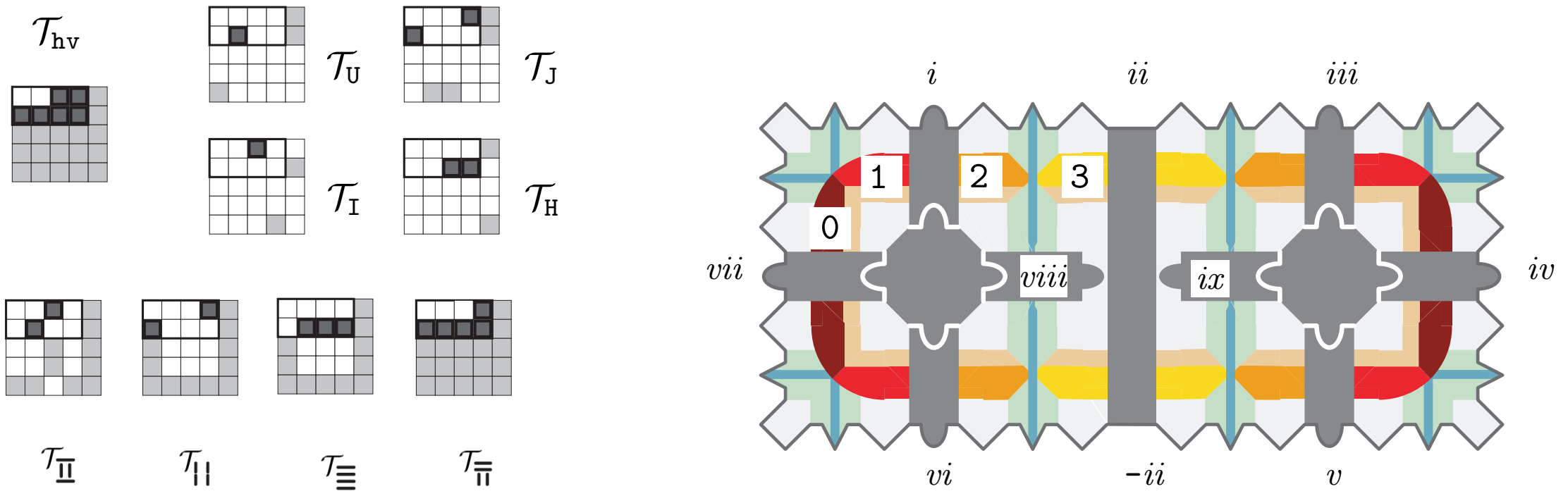

Figure 12: (Left) We indicate subsets of $\mathcal{T}_{\mathtt{hv}}$ — $\mathcal{T}_{\mathtt{H}}$, $\mathcal{T}_{\mathtt{I}}$, $\mathcal{T}_{\mathtt{J}}$ and $\mathcal{T}_{\mathtt{U}}$ — that are the building blocks of our aperiodic subsets of $\mathcal{T}_1$, four of which are shown at bottom: $\mathcal{T}_{\underline{\mathtt{II}}}$, $\mathcal{T}_{\mathtt{||}}$, $\mathcal{T}_{\underline{\equiv}}$ and $\mathcal{T}_{\overline{\mathtt{II}}}$.Note that $\mathcal{T}_{\mathtt{||}}$, $\mathcal{T}_{\underline{\equiv}}$, and $\mathcal{U}_1$ are all subsets of $\mathcal{T}_{\overline{\mathtt{II}}}$, but $\mathcal{T}_{\underline{\mathtt{II}}}$ and $\mathcal{T}_{\overline{\mathtt{II}}}$ each contain tiles the other does not. (Right) Indices for the markings on a marked block, the preimage of a domino under our natural local derivation, that we will use in substituting on marked blocks.

For each of the following, $\mathtt{X}(x)$ indicates a block of type $\mathtt{X}$, with $x$ being the marking $ii$. Working from Figure 12, for block

| | | |
|---|---|---|
| $\mathtt{U}(x)$ | we have | $i = iii = iv = vii = \mathtt{+}$; $v = vi = \mathtt{1}$; $viii = ix = \mathtt{0}$; $ii = x$; |
| $\mathtt{J}(x)$ | we have | $i = iii = v = vii = \mathtt{+}$; $iv = \mathtt{3}$; $vi = \mathtt{0}$; $viii = \mathtt{1}$; $ix = \mathtt{2}$; $ii = x$; |
| $\mathtt{I}(x)$ | we have | $i = iii = v = vi = \mathtt{+}$; $iv = vii = \mathtt{2}$; $viii = ix = \mathtt{3}$; $ii = x$; |
| $\mathtt{H}(x)$ | we have | $i = iii = \mathtt{2}$; $iv = vii = viii = ix = \mathtt{+}$; $v = vi = \mathtt{3}$; $ii = x$ . |

**The $\underline{\mathtt{II}}$ substitution** The substitution $\underline{\mathtt{II}}$ on $\mathtt{X}(x)$ gives

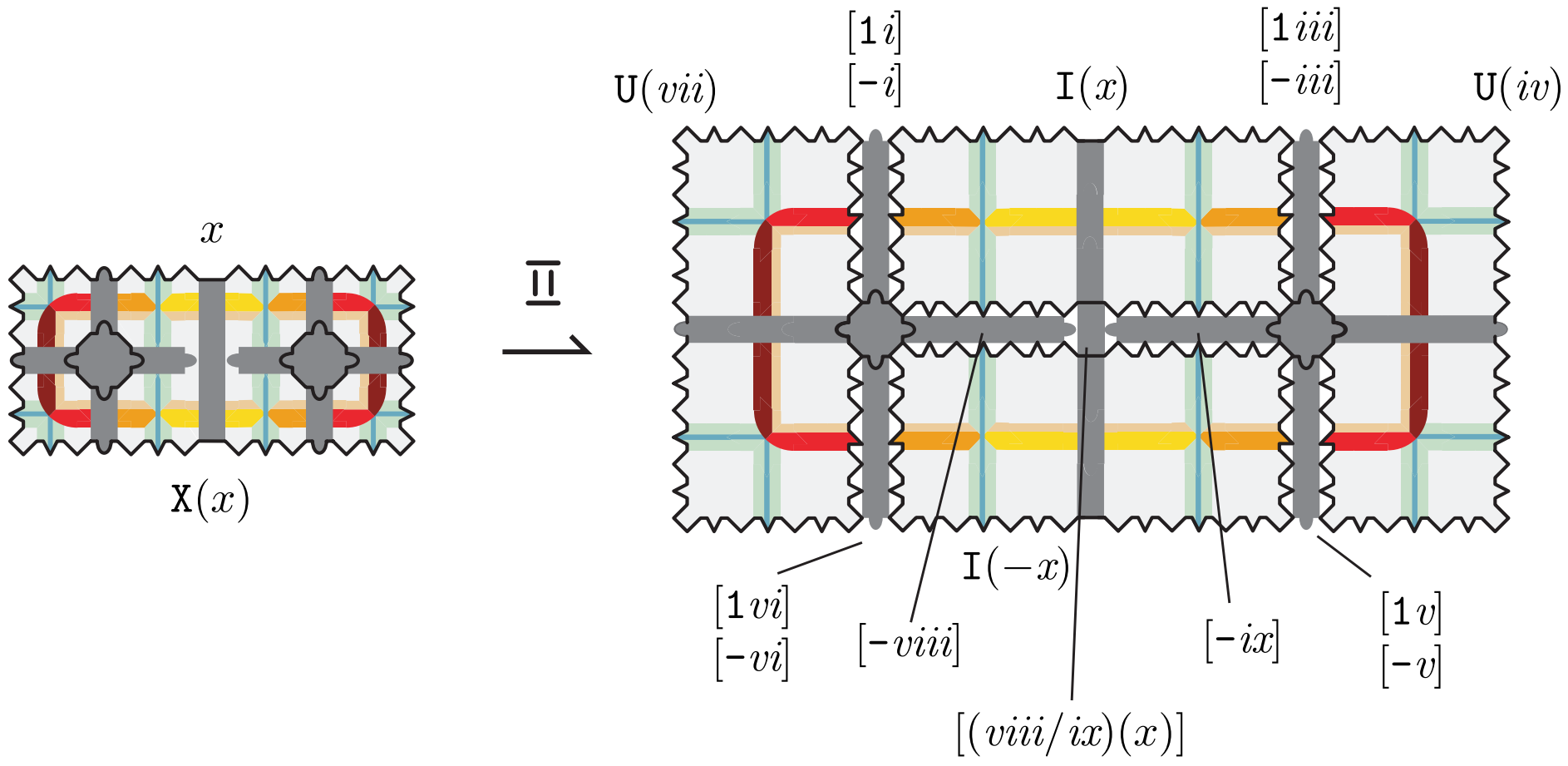

where $[(viii/ix)(x)]$ indicates $[viii\ x]$ or $[ix\ x]$ depending on the orientation and sidedness of the marking $x$, and the sidedness of the markings $viii$ and $ix$. The canonical notation for these tiles is easily resolved by examination.

Recalling that for $\mathtt{U}$, $i = iii = iv = vii = \mathtt{+}$; $v = vi = \mathtt{1}$; $viii = ix = \mathtt{0}$ and for $\mathtt{I}$, $i = iii = v = vi = \mathtt{+}$; $iv = vii = \mathtt{2}$; $viii = ix = \mathtt{3}$, we have

$$\mathtt{U}(x) \overset{\underline{\mathtt{II}}}{\mapsto} \mathtt{U}(\mathtt{+}) \quad \mathtt{I}(x) \quad \mathtt{I}(-x) \quad [\mathtt{11}] \quad [\mathtt{-0}] \quad [\mathtt{-1}] \quad [\mathtt{1+}] \quad [\mathtt{-+}] \quad [\mathtt{0}x]$$

$$\mathtt{I}(x) \overset{\underline{\mathtt{II}}}{\mapsto} \mathtt{U}(\mathtt{2}) \quad \mathtt{I}(x) \quad \mathtt{I}(-x) \quad [\mathtt{-3}] \quad [\mathtt{1+}] \quad [\mathtt{-+}] \quad [\mathtt{3}x] \quad [\mathtt{-3}x]$$

We iterate this further on U(+), U(2) and then I(+) and I(2). Obtaining no more, we halt fairly quickly, with just sixteen tiles, $\mathcal{T}_{\mathrm{II}}$. (Observe that [02] is the marked tile of Figure 8, arising precisely within a U(2) within any parent I($x$) block.)

$\mathcal{T}_{\mathrm{II}} := \mathcal{T}_{+} \cup \{$[02], [11], [32], [-32], [-0], [-1], [-3], [-+], [0+], [1+], [3+], [-3+]$\}$

**The ¦¦ substitution** Similarly, with the substitution ¦¦ on X($x$) we have

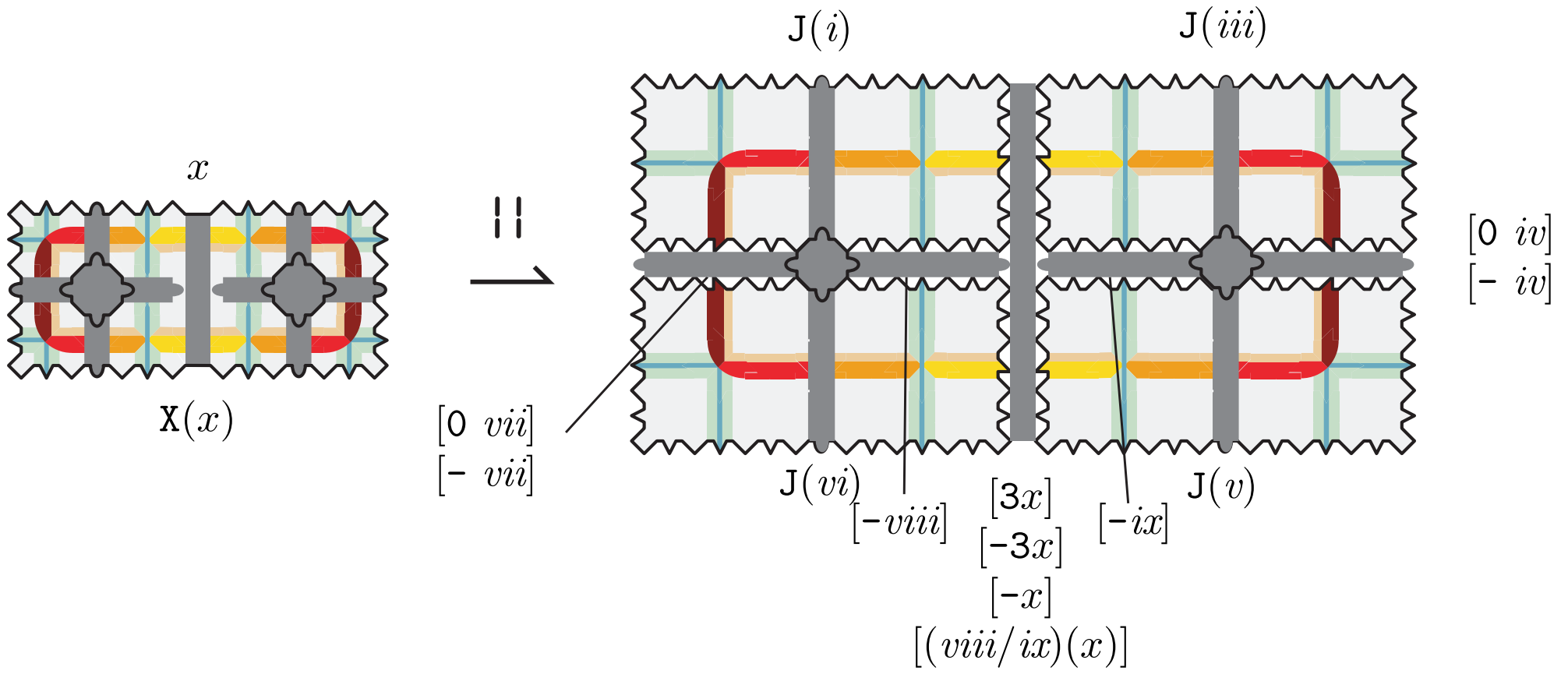

Noting that for J, $i = iii = v = vii =$ +; $iv =$ 3; $vi =$ 0; $viii =$ 1; $ix =$ 2 we have

J($x$) $\overset{¦¦}{\mapsto}$ J(0) J(+) [03] [0+] [-1] [-2] [-3] [-+] [1$x$] [3$x$] [-3$x$] [-$x$]

Further substituting on J(+) and J(0) we have additional tiles [10], [30], [-30], [-0], [1+], [3+] and [-3+]. We define:

$\mathcal{T}_{¦¦} := \mathcal{T}_{+} \cup \{$[03], [10], [30], [-30], [-0], [-1], [-2], [-3], [-+], [0+], [1+], [3+], [-3+]$\}$

**The ≣ substitution** For ≣, we have

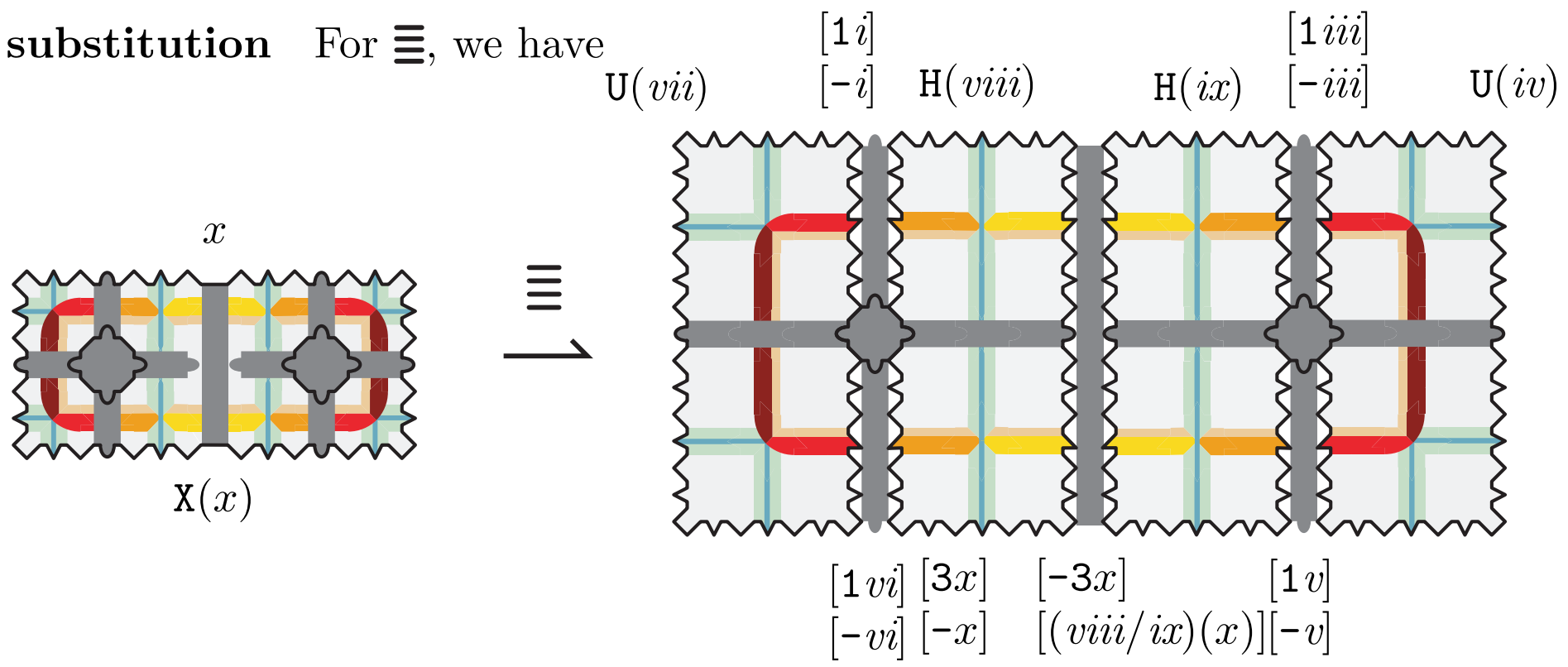

Since we have for U, $i = iii = iv = vii =$ +; $v = vi =$ 1; $viii = ix =$ 0 and for H, $i = iii =$ 2;

$iv = vii = viii = ix =$ +;   $v = vi =$ 3, we obtain

$$\mathtt{U}(x) \quad \overset{\equiv}{\mapsto} \quad \mathtt{U(+)} \quad \mathtt{H(0)} \quad [\mathtt{11}] \quad [\mathtt{-1}] \quad [\mathtt{1+}] \quad [\mathtt{-+}] \quad [\mathtt{0}x] \quad [\mathtt{3}x] \quad [\mathtt{-3}x] \quad [\mathtt{-}x]$$

$$\mathtt{H}(x) \quad \overset{\equiv}{\mapsto} \quad \mathtt{U(+)} \quad \mathtt{H(+)} \quad [\mathtt{12}] \quad [\mathtt{13}] \quad [\mathtt{-2}] \quad [\mathtt{-3}] \quad [\mathtt{3}x] \quad [\mathtt{-3}x] \quad [\mathtt{-}x]$$

Iterating once more on U(+), H(+) and H(0), we have tiles

$$\mathcal{T}_{\equiv} := \mathcal{T}_{+} \cup \{[\mathtt{11}], [\mathtt{12}], [\mathtt{13}], [\mathtt{30}], [\mathtt{-30}], [\mathtt{-0}], [\mathtt{-1}], [\mathtt{-2}], [\mathtt{-3}], [\mathtt{-+}], [\mathtt{0+}], [\mathtt{1+}], [\mathtt{3+}], [\mathtt{-3+}]\}$$

**The $\overline{\Pi}$ substitution** We face an additional complication when applying $\overline{\Pi}$: we have multiple cases, depending on the orientation of the substitution with respect to the tile; in particular, taking d to be one of the symmetries of the domino tile ⬚, d acts on the labels $i$, $ii$, ... $ix$.

On a generic block X($x$), we have

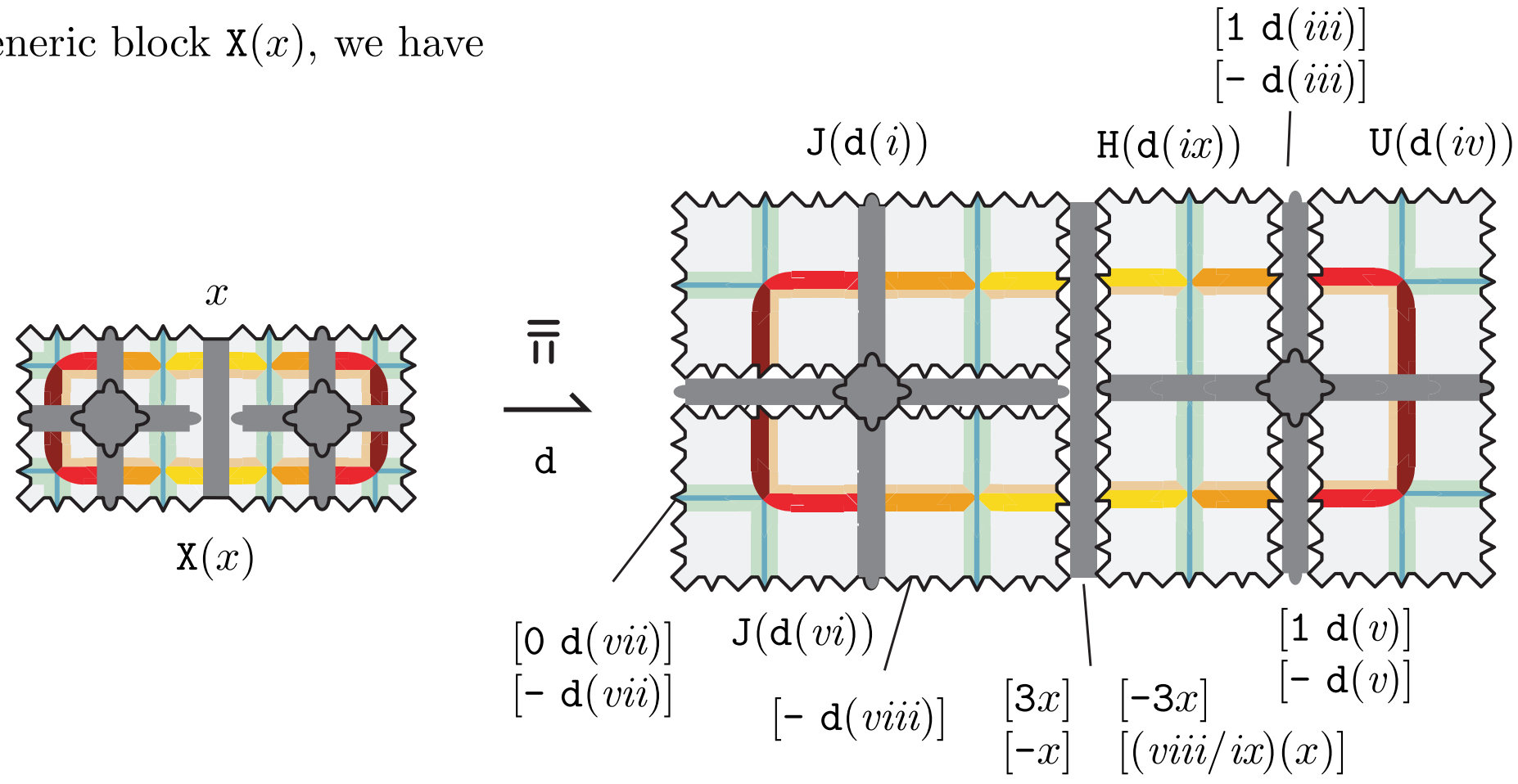

(Note the orientation of the marking $x$ must already be accounted for in reconciling tiles of the form [h$x$], regardless of the symmetry d.)

Because of the symmetries that H and U do have, irrespective of d,

$$\mathtt{H}(x) \quad \overset{\overline{\Pi}}{\mapsto} \quad \mathtt{J(2)} \quad \mathtt{J(3)} \quad \mathtt{H(+)} \quad \mathtt{U(+)} \quad [\mathtt{12}] \quad [\mathtt{13}] \quad [\mathtt{-2}] \quad [\mathtt{-3}] \quad [\mathtt{0+}] \quad [\mathtt{-+}] \quad [\mathtt{3}x] \quad [\mathtt{-3}x] \quad [\mathtt{-}x]$$

$$\mathtt{U}(x) \quad \overset{\overline{\Pi}}{\mapsto} \quad \mathtt{J(1)} \quad \mathtt{J(+)} \quad \mathtt{H(0)} \quad \mathtt{U(+)} \quad [\mathtt{11}] \quad [\mathtt{-1}] \quad [\mathtt{0+}] \quad [\mathtt{1+}] \quad [\mathtt{-+}] \quad [\mathtt{0}x] \quad [\mathtt{3}x] \quad [\mathtt{-3}x] \quad [\mathtt{-}x]$$

When d preserves left and right in the figure above (that is, d = 0, 2 in the notation of the following section) we have

$$\mathtt{J}(x) \quad \overset{\overline{\Pi}}{\mapsto} \quad \mathtt{J(0)} \quad \mathtt{J(+)} \quad \mathtt{H(2)} \quad \mathtt{U(3)} \quad [\mathtt{-1}] \quad [\mathtt{0+}] \quad [\mathtt{1+}] \quad [\mathtt{-+}] \quad [\mathtt{1}x] \quad [\mathtt{3}x] \quad [\mathtt{-3}x] \quad [\mathtt{-}x]$$

and we further iterate on J(+, 0, 1, 2, 3), H(+, 0, 2), U(+, 3).

When d reverses left and right in the figure above (that is, d = 1, 3) we have

$$\mathtt{J}(x) \quad \overset{\overline{\Pi}}{\mapsto} \quad \mathtt{J(+)} \quad \mathtt{H(1)} \quad \mathtt{U(+)} \quad [\mathtt{03}] \quad [\mathtt{10}] \quad [\mathtt{-0}] \quad [\mathtt{-2}] \quad [\mathtt{-3}] \quad [\mathtt{1+}] \quad [\mathtt{-+}] \quad [\mathtt{1}x] \quad [\mathtt{3}x] \quad [\mathtt{-3}x] \quad [\mathtt{-}x]$$

and we iterate on J(+, 1, 2, 3), H(+, 0, 1), U(+).

However in either case, regardless of d, we require all of the following tiles, and no more:

$\mathcal{T}_{╥} := \mathcal{T}_1 \setminus \{[02]\}$

## 3.5 Concluding the proof of Theorem 1

We observe that the substitution

Ⅱ includes blocks I, U, but not blocks J, H;

¦¦ includes blocks J, but not blocks H, I, U;

≣ includes blocks H, U, but not I, J; and

╥ includes blocks H, J, U, but not I.

We will say that sets $A$ and $B$ are "mutually excludable" if $A \setminus B$ and $B \setminus A$ are both non-empty, or equivalently that $A \not\subset B$ and $B \not\subset A$. The set of tiles

$\mathcal{T}_{Ⅱ}$ contains $\mathcal{T}_{\mathtt{I}}$ and $\mathcal{T}_{\mathtt{U}}$ but is mutually excludable with $\mathcal{T}_{\mathtt{H}}$ and $\mathcal{T}_{\mathtt{J}}$, hence with $\mathcal{T}_{¦¦}$, $\mathcal{T}_{≣}$, $\mathcal{T}_{╥}$.

$\mathcal{T}_{¦¦}$ contains $\mathcal{T}_{\mathtt{J}}$ but is mutually excludable with $\mathcal{T}_{\mathtt{H}}$, $\mathcal{T}_{\mathtt{I}}$ and $\mathcal{T}_{\mathtt{U}}$, hence with $\mathcal{T}_{Ⅱ}$ and $\mathcal{T}_{≣}$.

$\mathcal{T}_{≣}$ contains $\mathcal{T}_{\mathtt{H}}$ and $\mathcal{T}_{\mathtt{U}}$ but is mutually excludable with $\mathcal{T}_{\mathtt{I}}$ and $\mathcal{T}_{\mathtt{J}}$, hence with $\mathcal{T}_{Ⅱ}$ and $\mathcal{T}_{¦¦}$.

$\mathcal{T}_{╥}$ contains $\mathcal{T}_{\mathtt{H}}$, $\mathcal{T}_{\mathtt{J}}$ and $\mathcal{T}_{\mathtt{U}}$ but is mutually excludable with $\mathcal{T}_{\mathtt{I}}$ hence with $\mathcal{T}_{Ⅱ}$.

On the other hand, $\mathcal{T}_{¦¦}$, $\mathcal{T}_{≣}$, and $\mathcal{U}_1$ are each subsets of $\mathcal{T}_{╥}$.

Consequently

| the set of tiles | admits blocks | but does not admit | hence admits supertiles | but is insufficient to admit supertiles | and thus can and only can enforce substitution |
|---|---|---|---|---|---|
| $\mathcal{T}_{Ⅱ}$ | I, U | H, J | Ⅱ | ¦¦, ≣, ╥ | Ⅱ |
| $\mathcal{T}_{¦¦}$ | J | H, I, U | ¦¦ | Ⅱ, ≣, ╥ | ¦¦. |
| $\mathcal{T}_{≣}$ | H, U | I, J | ≣ | Ⅱ, ¦¦, ╥ | ≣. |
| $\mathcal{T}_{╥}$ | H, J, U | I | ¦¦, ≣, ╥ | Ⅱ | {¦¦, ≣, ╥} . |

Taking unions of these sets,

| the set of tiles | admits blocks | but does not admit | hence admits supertiles | but is insufficient to admit supertiles | and thus can and only can enforce substitution |
|---|---|---|---|---|---|
| $\mathcal{T}_{Ⅱ} \cup \mathcal{T}_{¦¦}$ | I, J, U | H | Ⅱ, ¦¦ | ≣, ╥ | {Ⅱ, ¦¦} |
| $\mathcal{T}_{Ⅱ} \cup \mathcal{T}_{≣}$ | H, I, U | J | Ⅱ, ≣ | ¦¦, ╥ | {Ⅱ, ≣} |
| $\mathcal{T}_{¦¦} \cup \mathcal{T}_{≣}$ | H, J, U | I | ¦¦, ≣ | Ⅱ, ╥(*) | {¦¦, ≣} |
| $\mathcal{T}_{Ⅱ} \cup \mathcal{T}_{¦¦} \cup \mathcal{T}_{≣}$ | H, I, J, U | | Ⅱ, ¦¦, ≣ | ╥(*) | {Ⅱ, ¦¦, ≣} |
| $\mathcal{T}_1$ | H, I, J, U | | Ⅱ, ¦¦, ≣, ╥ | | {Ⅱ, ¦¦, ≣, ╥}. |

Just on the basis of which blocks each set of tiles can admit, we will have thus established nearly every case stated in the theorem. The cases marked (*) above require one further observation: $\mathcal{T}_{¦¦} \cup \mathcal{T}_{≣}$ and $\mathcal{T}_{Ⅱ} \cup \mathcal{T}_{¦¦} \cup \mathcal{T}_{≣}$ are *proper* subsets of $\mathcal{T}_{╥}$ and thus cannot admit ╥ blocks.

Iterating the construction of larger and larger blocks of the types admitted by these subsets, we can form arbitrarily large configurations that are locally derivable to larger and larger supertiles in the corresponding tiling substitution system.

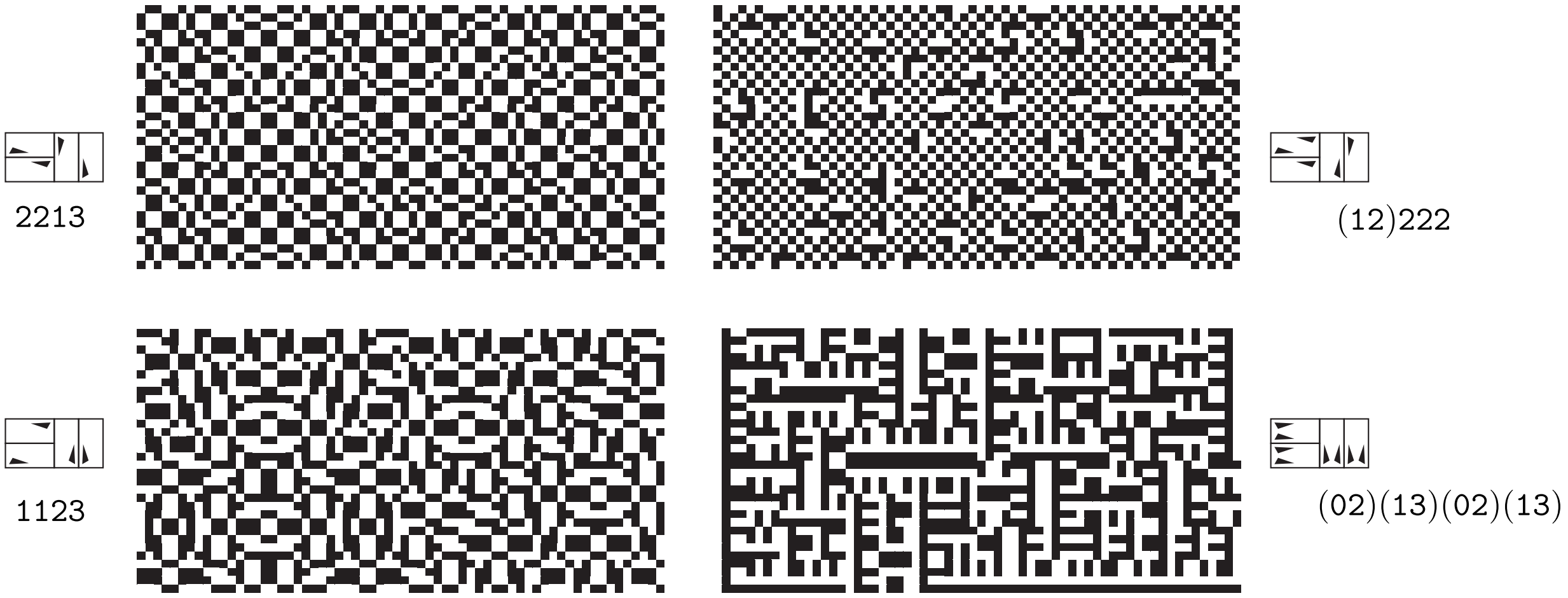

Figure 13: A few examples of $\overline{\Pi}$ tiling substitution systems, with dominoes colored as ▮ ▭. Examining a fuller catalogue, as in [6], one notices a wide variety of qualitative differences in statistical structure, which awaits quantification.

Each of these sets *can* tile the plane, because it admits configurations covering arbitrarily large disks — this standard argument appears as The Extension Theorem, Theorem 3.8.1 in [8] and is often cited as "by compactness" or "by Koenig's Lemma".

On the other hand, these tiles admit no other tilings whatsoever.

Consider any subset $\mathcal{T}$ of $\mathcal{T}_1$. We have shown that in order to admit a tiling, $\mathcal{T}$ must include all of the outward tiles of $\mathcal{T}_+$, and that in any tiling by $\mathcal{T}$, the cornered outward tiles can only lie in blocks `H`, `I`, `J` and `U`; these can only be admitted if the corresponding sets $\mathcal{T}_{\mathtt{H}}, \mathcal{T}_{\mathtt{I}}, \mathcal{T}_{\mathtt{J}}$, or $\mathcal{T}_{\mathtt{U}}$ are included in $\mathcal{T}$. These must lie in larger and larger marked supertiles, *ad infinitum*. The only marked supertiles that any of these sets admits are those indicated in the table above. That is, every tiling in $\Sigma(\mathcal{T}_{\overline{\underline{\Pi}}})$ is locally decomposible to a hierarchical tiling defined by the corresponding tiling substitution system.

Moreover, within these tilings, these hierarchies are unique, and the tilings are thus non-periodic. These tile sets are *aperiodic*

Finally, a few last observations: Just as described in the proof of the first case, in each of our named aperiodic sets $\mathcal{T} \subset \mathcal{T}_1$, $\mathcal{T} \cup \mathcal{U}_1$ admits just the same tilings, with any additional tiles in $\mathcal{U}_1 \setminus \mathcal{T}$ unused and wasted. On the other hand, in order to admit the entire such space, and not restrict which rules may be applied, every tile in $\mathcal{T}$ appears after iterating the substitutions, and so is essential. And in particular, no proper subset of $\mathcal{T}_{\overline{\text{II}}}$, $\mathcal{T}_{\text{¦¦}}$ or $\mathcal{T}_{\equiv}$ even admits any tiling at all.

Though the ¦¦ (and ≣) substitution tilings are periodic, the matching rule tilings enforce heirarchies — each corner tile lies in a larger and larger, unique, hierarchy of marked supertiles. The tilings the the tiles admit are non-periodic and so the sets of tiles in Theorem 1 are aperiodic.

## 4 25380 aperiodic subsets of 211 tiles

Theorem 2 produces, for any of the $\overline{\Pi}$ substitution systems, a set of tiles that enforce it, as illustrated in Example 4.5. Our goal in this section is to define our terms in order to state and

prove:

**Theorem 2** *For each $\overline{\Pi}$ tiling substitution $\sigma_{\mathtt{S}}$, with symbol $\mathtt{S}$, the set of tiles*

$$\mathcal{T}_{\mathtt{S}} := \mathcal{T}_{+_2} \cup \mathcal{T}_0 \cup \left( \bigcup_{\alpha \in \mathtt{S}} \mathcal{T}_\alpha \right) \cup \left( \bigcup_{\alpha,\beta \in \mathtt{S}} \mathcal{T}_{\beta\to\alpha} \right)$$

*is aperiodic, enforcing the substitution $\sigma_{\mathtt{S}}$. Moreover, any subset of $\mathcal{U}_2 \cup \mathcal{T}_{\mathtt{S}}$ that enforces $\sigma_{\mathtt{S}}$ must contain $\mathcal{T}_{\mathtt{S}}$. If $\mathtt{S}$ is deterministic, no proper subset of $\mathcal{T}_{\mathtt{S}}$ admits any tiling.*

Essentially Theorem 2 states that any given a $\overline{\Pi}$ substitution $\sigma_{\mathtt{S}}$, denoted by a to-be-defined symbol $\mathtt{S}$, is precisely enforced by the union of the "atomic sets of tiles" corresponding to the "atomic symbols" making up $\mathtt{S}$, together with some standard cornered and crossing tiles.

One may safely have passed over the latter half of Section 3; here we use only that the set of tiles $\mathcal{T}_{\overline{\Pi}}$ admits exactly the blocks `J`, `H` and `U` of Figure 11.

Section 4.1 describes the `d`-channel markings `0`, `1`, `2` and `3`. Section 4.2 presents a system for encoding $\overline{\Pi}$ substitution rules and the atomic substitutions of which they are comprised. In Section 4.3 we define our master set $\mathcal{T}_2$ and its subsets $\mathcal{T}_0, \mathcal{T}_{+_2}$. In Section 4.4 we define the atomic sets of tiles, the sets $\mathcal{T}_\alpha$, $\mathcal{T}_{\beta\to\alpha}$ of the theorem, with fully worked out examples in Section 4.5.

Finally, in Section 4.6, the theorem now being well-defined, we complete its proof, which has essentially the same outline as in Section 3.1, that of [2] onward, although we must fully industrialize the arguments: For each $\overline{\Pi}$ substitution system we show that the defined tiles can form, and are necessary to form, arbitrarily large marked supertiles. In Section 5 we carefully count these sets up to symmetry.

## 4.1 The `d`-channel markings

As illustrated in Figure 14 and defined in Section 2, our `d`-channel markings are $\mathtt{0}, \mathtt{1}, \mathtt{2}, \mathtt{3}$, in $\mathbb{Z}_2 \oplus \mathbb{Z}_2$ under nim addition, that is, addition, without carry, on binary strings (so for example $\mathtt{2} + \mathtt{3} = (\mathtt{10}) + (\mathtt{11}) = (\mathtt{01}) = \mathtt{1}$). These correspond to the symmetries of the domino tile, within some local framing: nim adding `1` corresponds to a reflection across a vertical mirror line; nim adding `2` to a reflection across a horizontal one; taking both together, nim adding `3` corresponds to a two-fold rotation; and of course nim adding `0` accomplishes nothing.

As sketched in Figure 15, with these new markings we can first control where the different kinds of blocks within a supertile must lie: in Section 4.3 we will define "key tiles", special crossing tiles $\mathcal{K}_2$, lying at the circled locations in the figure, where signed markings meet. These key tiles exactly permit a pair of `J` blocks in the `0` and `1` quarters of a supertile (by using only the `0` and `1` markings in the `d`-channel of the vertical markings of the key tiles of $\mathcal{T}_{\mathtt{J}}$) and similarly allow a pair of $\mathtt{H}, \mathtt{U}$ blocks on the `2` and `3` quarters of a supertile (by using only the `2` and `3` markings in the `d`-channel of the vertical markings of the key tiles of $\mathcal{T}_{\mathtt{H}}$ and of $\mathcal{T}_{\mathtt{U}}$). We arbitrarily take this arrangement as canonical.

There are various possibilities for the horizontal `d`-channel markings for a tile in $\mathcal{K}_2$, depending on its location in the marked $\overline{\Pi}$ configuration of Figure 15. For example the pair tiles at left both appear at a short end of a `J` block, and so their horizontal markings might be, and could only be, a `0`-`1` pair or a `2`-`3` pair. Similarly, the eight tiles at right might have, and could only have, `0` or `1` or `2` or `3` for the `d`-channel of the horizontal marking. In Figure 18 the key tiles $\mathcal{K}_2$ are illustrated within $\mathcal{T}_2$; and they are defined more precisely on page 24.

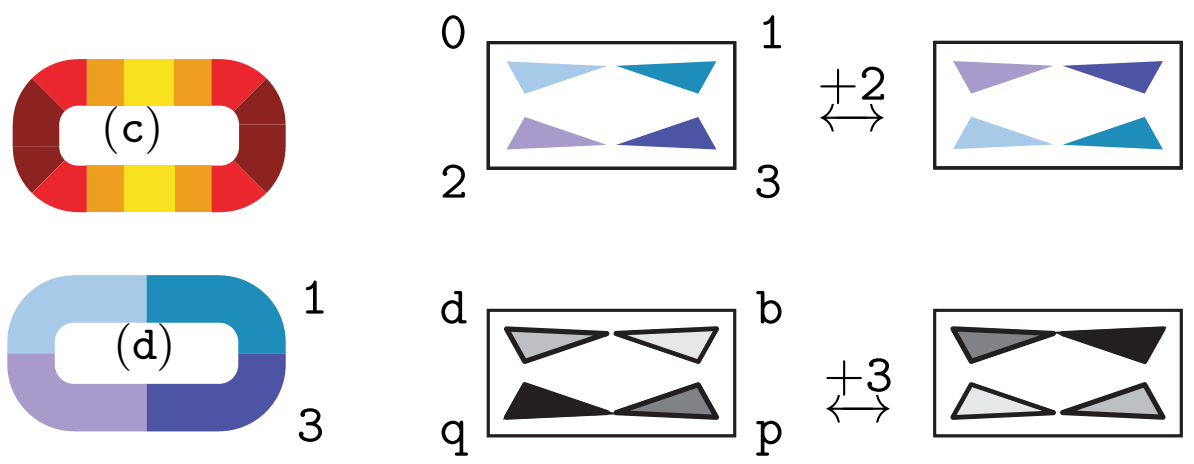

| + 0 | +1 | +2 | +3 |
|---|---|---|---|
| 1 | 0 | 3 | 2 |
| 2 | 3 | 0 | 1 |
| 3 | 2 | 1 | 0 |
| d | b | q | p |

Figure 14: The d-channel markings 0, 1, 2 and 3 are in $\mathbb{Z}_2 \oplus \mathbb{Z}_2$ under nim addition, framing the orientation of each supertile. The markings naturally act upon each other and upon the variables d, b, p, q for generic framings: For example, at top center, we reflect a domino top to bottom, nim adding +2 to the markings 0, 1, 2 and 3, the results of which we can see in the figure or read off in the +2 column (or row) of the table at right. At bottom center, +3 acts on d, b, q, p and on the domino by a two-fold rotation. As these digits act on the d-channel, they naturally act on markings, tiles, tilings and tiling spaces, as we use in defining our tiles in §4.3 and 4.4, and in counting tiling matching rule spaces up to symmetry in Section 5.

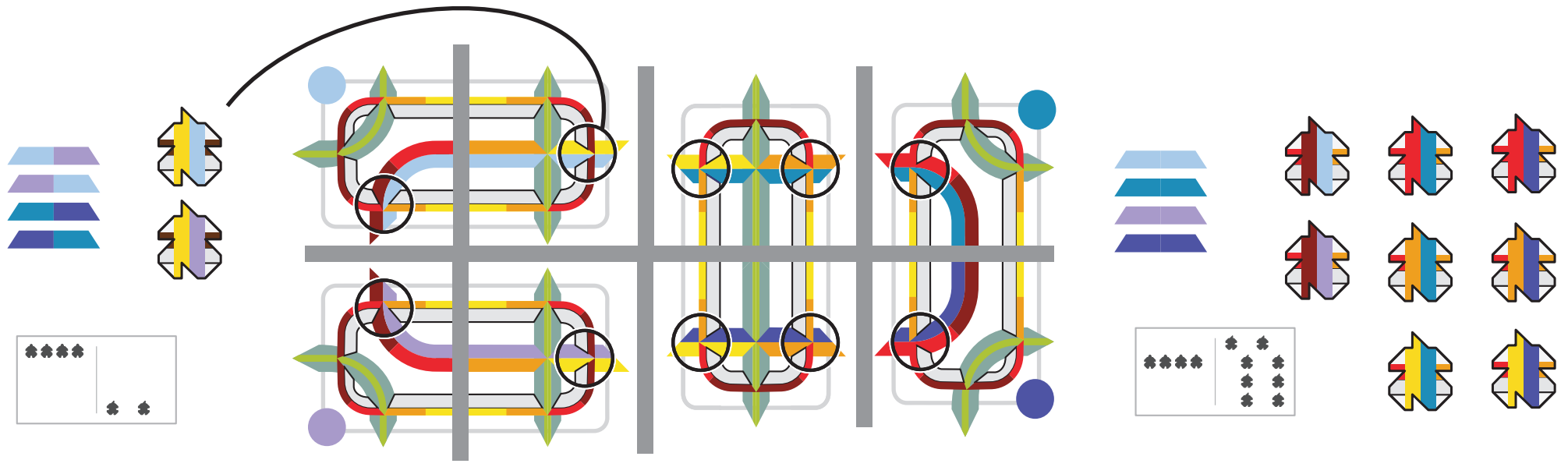

Figure 15: We (arbitrarily but canonically) frame each ╦ block with d-channel markings as shown, with a J block marked 0, another marked 2, and a H, U pair together marked 1 and 3. The "key tiles" lie at the circled locations, where signed markings meet. We restrict the allowed vertical (cd) and horizontal (c) markings, allowing a few possibilities for the horizontal (d), those shown, so that we can inductively force our supertiles to be well-formed ╦ blocks as we discuss further in §4.6. We may further restrict the orientations of the children with respect to that of the parent, by further limiting *which* of these horizontal d-channel markings to permit— these are our atomic subsets of $\mathcal{K}_2$. The small diagrams at bottom indicate the locations of the key tiles in the table of Figure 18, which shows all of $\mathcal{T}_2$.

By further restricting *which* horizontal d-channel markings to include or exclude in a set of tiles from $\mathcal{K}_2$, we can determine how these children blocks can be oriented within a parent ╦ configuration.

To this end, we further use these digits 0, 1, 2, 3 to define an *ad hoc* (but intentionally chosen) encoding of how a child supertile may fit within a parent supertile, that is, the relative orientations of each with respect to each other. For each of the four children, we arbitrarily fix a specific orientation to be 0, and define the other orientations for that child accordingly, as indicated in Figure 16.

In a ╦ supertile, for each of the four children, denoted s, t, u and v, we specify which orientations we will allow, defining a code for the the ╦ substitution tiling systems, described shortly in Section 4.2. Such a code S will be composed of elementary, "atomic" symbols $\alpha, \beta \in$ S that each describe a portion of what is allowed in the full substitution rule $\sigma_{\mathtt{S}}$ given by S.

For example, the symbol ·320 (itself only a "partial" ╦ substitution rule) is shown at left in Figure 17; it is composed of atomic symbols ·3··, ··2·, and ···0.

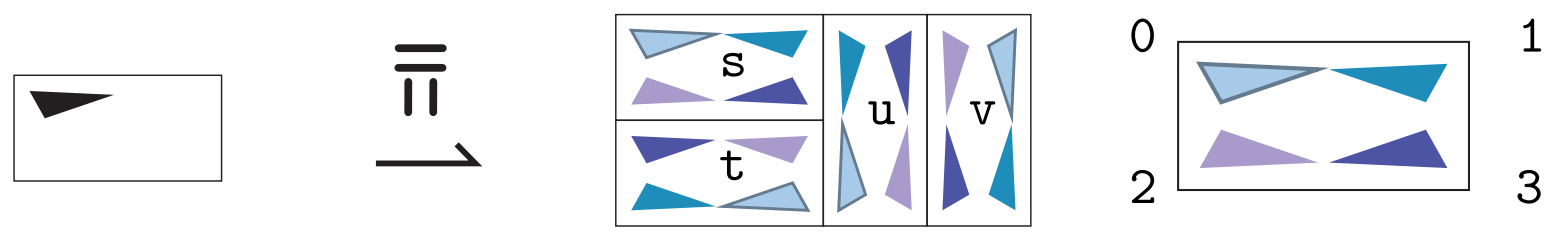

Figure 16: The d-channel markings encode the relative orientations of the children supertiles with respect to a parent supertile, allowing us to denote any $\overline{\overline{\Pi}}$ substitution rule as a quadruple $\mathtt{S} = \mathtt{stuv}$ of subsets of our digits: At right in Figure 3 we show 1101; the partial rule (with $\mathtt{s} = \emptyset$) is shown at left in Figure 17; and a few more examples appear in Figure 13. Catalogues of examples appear in [6] and an interactive *Mathematica* demonstration is at [7].

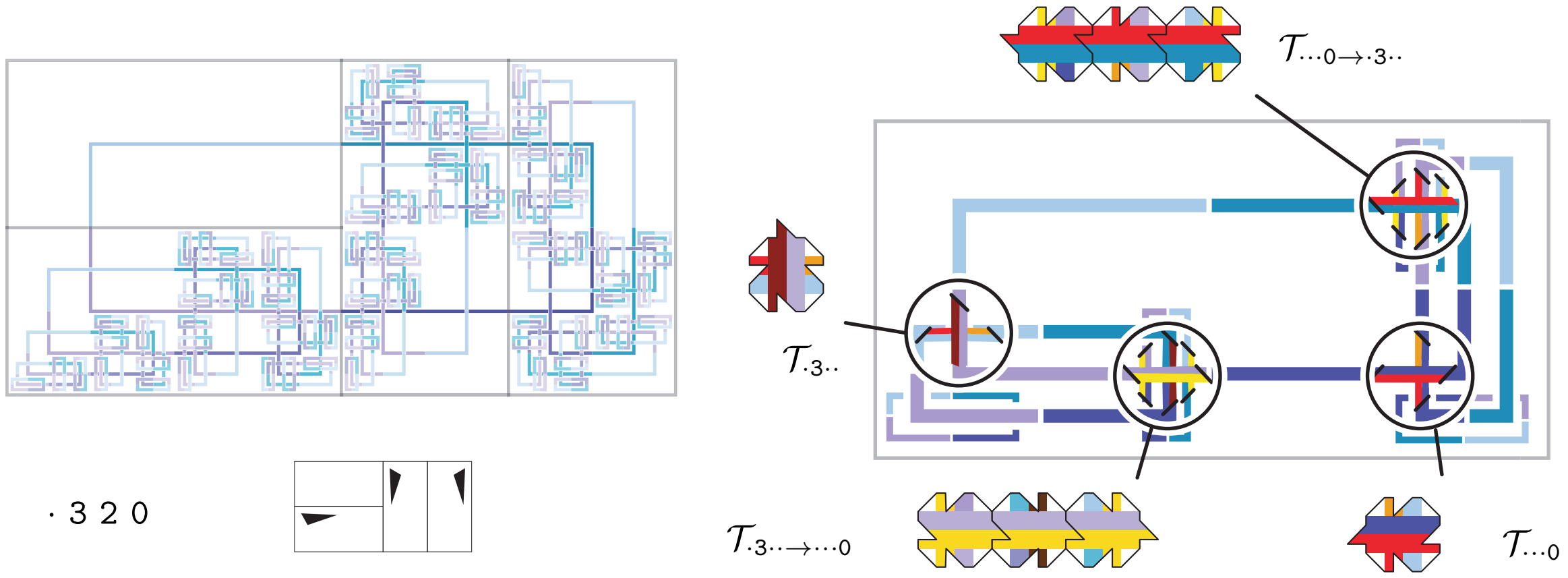

Figure 17: In the compressed notation described in Section 4.2, the system · 3 2 0 is shown at left. At right, the atomic subsets necessary and sufficient (with a full complement of other atomic sets) to include from $\mathcal{T}_2$ in order to enforce ·3·0.

In Section 4.4, for any atomic symbols $\alpha, \beta$, we define a corresponding atomic set of tiles $\mathcal{T}_\alpha, \mathcal{T}_{\beta\to\alpha}$ each of which is necessary and sufficient to allow a particular relationship between a child and parent supertile. At right in Figure 17, we see the atomic sets needed for any substitution that includes the partial symbol ·3·0, namely $\mathcal{T}_{\cdot 3\cdot\cdot}$, $\mathcal{T}_{\cdot\cdot\cdot 0}$, $\mathcal{T}_{\cdot 3\cdot\cdot\to\cdot\cdot\cdot 0}$, and $\mathcal{T}_{\cdot\cdot\cdot 0\to\cdot 3\cdot\cdot}$.

Once we define our terms, further below, Theorem 2 states that a given a $\overline{\overline{\Pi}}$ substitution $\sigma_{\mathtt{S}}$ with symbol S is precisely enforced by the union of the atomic sets of tiles corresponding to the atomic symbols making up S, together with some standard cornered and crossing tiles.

An example of atomic sets of tiles is illustrated at upper right in Figure 17. In Example 4.5 we calculate $\mathcal{T}_{\mathtt{S}}$ for S = 1101, with 67 tiles, and 1023, with 65.

Catalogues of colour figures of $\overline{\overline{\Pi}}$ substitution tilings are in [6], and the *Mathematica* demonstration at [7] is useful for gaining intuition. Though we actually derive 50625 aperiodic sets below, most of these occur in MLD pairs, which we carefully out count in Section 5.

### 4.2 Symbols for $\overline{\overline{\Pi}}$ tiling substitution rules:

We define symbols describing $\overline{\overline{\Pi}}$ tiling substitution rules, and their constituent, atomic, partial tiling substitution rules. (Again, see [6] for many examples).

Each symbol S is a list of four subsets s, t, u, v of $\{0, 1, 2, 3\}$, including possibly the empty set; each of these subsets gives the allowed orientations of a corresponding child within a parent supertile, as shown in Figure 16. We take many notational conveniences. Any singleton as $\{a\}$ is denoted more simply as $a$. Any empty digit $\{\}$ we denote $\cdot$. If every digit $\{\mathtt{0},\mathtt{1},\mathtt{2},\mathtt{3}\}$ is used, we write $*$. Any other subset $\{a,\dots,z\}$ of digits, we write as $(a..z)$. The **full** substitutions are $\overline{\sqcap\!\sqcap}$ tiling substitutions with no empty digits, and are precisely those that give well-defined $\overline{\sqcap\!\sqcap}$ substitutions in the sense of Section 1. The full substitutions with all digits singletons are **deterministic**.

For example, if $\mathtt{S} = (\mathtt{12})\mathtt{222}$, as at top right of Figure 13, then $\mathtt{s} = \{\mathtt{1},\mathtt{2}\}$, $\mathtt{t} = \{\mathtt{2}\}$, $\mathtt{u} = \{\mathtt{2}\}$ and $\mathtt{v} = \{\mathtt{2}\}$. As formally defined as in Section 1 there two deterministic tiling substitution rules, denoted 1222, 2222 which together give a (non-deterministic) tiling substitution system $\sigma_{(\mathtt{01})\mathtt{231}}$. (The system at bottom right of Figure 13 has $2^4$ constituent deterministic rules.)

The **atomic** substitutions are partial substitutions with just one singleton among its digits with all the other digits being empty. Given a symbol S, we will write that the atomic substitution $\alpha \in \mathtt{S}$ if its one non-empty digit is one of the digits in S. For example, the atomic substitutions $\mathtt{0}\cdots$, $\mathtt{1}\cdots$, $\cdot\mathtt{2}\cdot\cdot$, $\cdot\cdot\mathtt{3}\cdot$ and $\cdots\mathtt{1}$ are exactly those that are in the full substitution (01)231.

The symmetries of the domino naturally operate on symbols in the d-channel, and in Section 5 we find useful notation such as m+2 meaning $(\mathtt{abc}(\mathtt{d}{+}\mathtt{2}))_{\mathtt{m}}$ (and so too, for example, by definition (+) + 2 = (+)).

### 4.3 The second set of tiles $\mathcal{T}_2$ and its basic subsets $\mathcal{T}_{+_2}$, $\mathcal{K}_2$, $\mathcal{U}_2$, $\mathcal{T}_0$

In Figure 18 we give notation for our markings and our tiles in $\mathcal{T}_2$:

We define sixteen outward tiles $\mathcal{T}_{+_2}$, exactly those of $\mathcal{T}_{+}$, recalling the notation of Section 2, with a choice for the value of $\mathtt{d} = \mathtt{0},\mathtt{1},\mathtt{2},\mathtt{3}$ — that is, our outward tiles maybe cornered or uncornered and have markings
(+)(+)(++1d)(+-0d) or (+)(+-2d)(+)(++3d).

We name our new cross tiles $[xy|zw]$, in
[ -, (013(-3))(0123) | +, (0123)(0123) ] — $x$ and $z$ are the abbreviations for the horizontal and vertical markings of Figure 9, but we are no longer suppressing the marking d, which now may be any of $\mathtt{0},\mathtt{1},\mathtt{2},\mathtt{3} =: y, w$.

As before if the horizontal marking is the plain (-000), we write - for $xy$; otherwise:

| west | and east | denoted |
|---|---|---|
| (--0d) | (-+0q) | [0d\| |
| (--1d) | (-+2d) | [1d\| |
| (--3d) | (-+3b) | [3d\| |
| (-+3d) | (--3b) | [-3d\| |

Similarly, we denote north and south vertical markings (++$zw$) and (--$zw$) as $|zw]$. In the special case that the north marking is ($\pm$000), we take $|zw]$ as |+].

Together these give names for our cross tiles of the form $[xy|zw]$. These names are unique with a few exceptions: In some cases, such as [30|+] = [31|+], we take the lexigraphically earlier name to be canonical.

It will be useful to apply digits in $\mathtt{0},\mathtt{1}\dots,\mathtt{d}$, etc. to our markings and so to our tiles in $\mathcal{T}_2$, writing, $[xy|zw] + \mathtt{d}$ to mean $[x(y+\mathtt{d})|z(w+\mathtt{d})]$, and thus to our sets of tiles and the tiling spaces they enforce.

For example, we will write that $\mathcal{T}_{\cdot\cdot\mathtt{d}\cdot}$ includes [1b|31], meaning that, taking d as 0, 1, 2, 3:

$\mathcal{T}_{\cdot\cdot\mathtt{0}\cdot}$ includes [11|31];

$\mathcal{T}_{\cdot\cdot\mathtt{1}\cdot}$ includes [10|31];

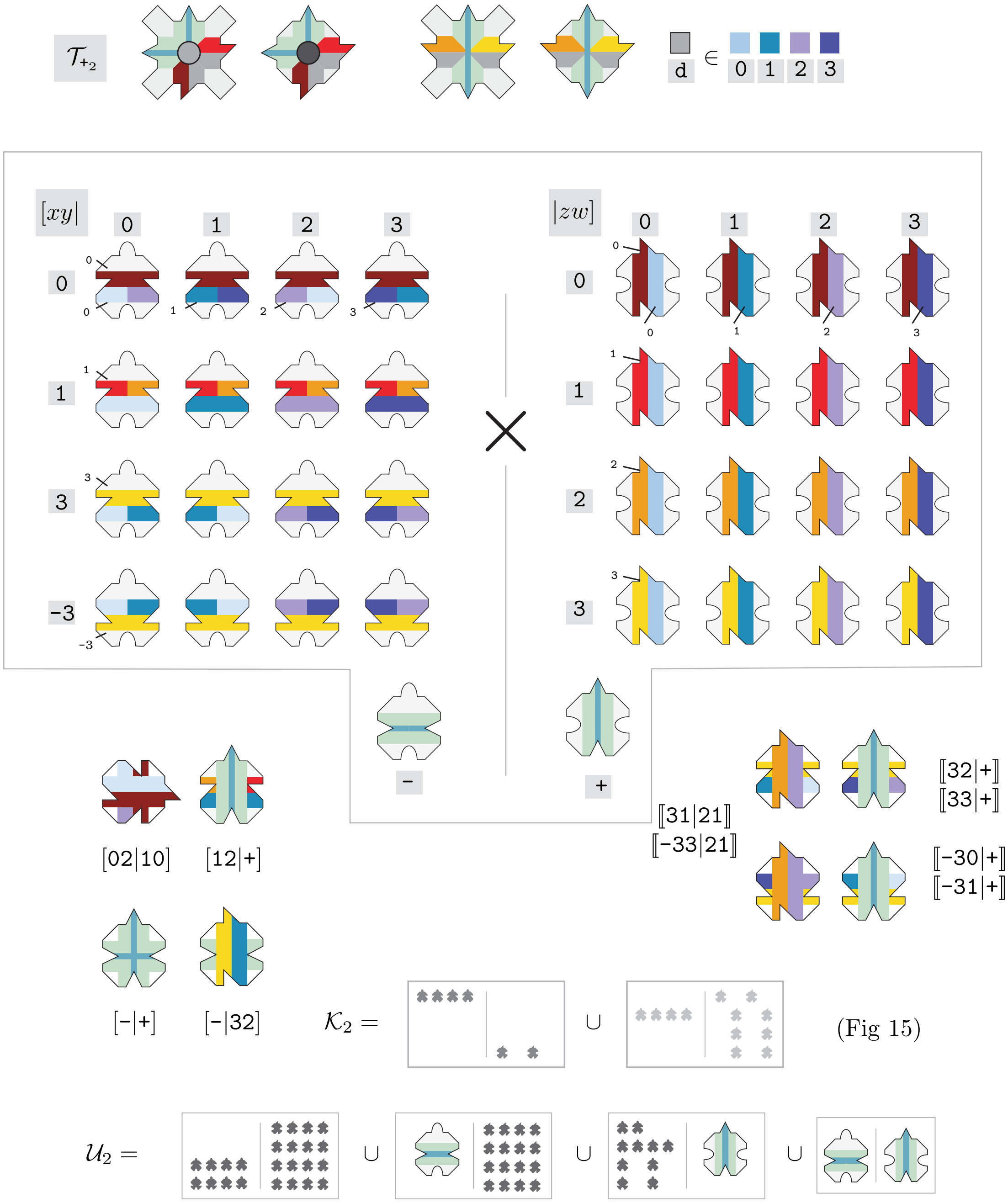

Figure 18: Our second set of tiles $\mathcal{T}_2$, the tiles of Figure 9), but using `d` = `0, 1, 2, 3` in our markings. At top sixteen outward tiles, $\mathcal{T}_{+_2}$. The cross tiles, indicated in the table, are named [$xy|zw$], where [$xy$| denotes the horizontal markings, and |$zw$] the vertical ones. (We give a few for practice.) Some names describe the same tile — for example [00|+] = [02|+] and we always use the lexigraphically earlier one. Tiles [3d|$zw$] and [-3q|$zw$] will always be used in pairs, ⟦3d|$zw$⟧ = ⟦-3q|$zw$⟧, such as those named here. We use only 195 of the possible cross tiles in $\mathcal{T}_2$: At bottom the sets of key tiles $\mathcal{K}_2$ — which particularly determine which systems are enforced — and the rest, $\mathcal{U}_2$.

$\mathcal{T}_{\cdot\cdot 2\cdot}$ includes $[\texttt{13}|\texttt{31}]$; and
$\mathcal{T}_{\cdot\cdot 3\cdot}$ includes $[\texttt{12}|\texttt{31}]$.

Tiles with $x = \texttt{3}$ or $\texttt{-3}$ will always be used in matching pairs, which we denote

$$[\![\texttt{3d}|zw]\!] := [\![\texttt{-3q}|zw]\!] := \{[\texttt{3d}|zw], [\texttt{-3q}|zw]\}$$

taking the first notation as canonical. Note that $[\![\texttt{30}|\texttt{+}]\!] = [\![\texttt{31}|\texttt{+}]\!]$ and $[\![\texttt{32}|\texttt{+}]\!] = [\![\texttt{33}|\texttt{+}]\!]$, and as always, we take the lexigraphically earlier name.

We define forty key tiles, indicated at the center upper bottom of Figure 18, and described and motivated in Figure 15.

$\mathcal{K}_2 :=$ $\{[\texttt{0d}|\texttt{30}], [\texttt{0d}|\texttt{32}], [\texttt{1d}|\texttt{00}], [\texttt{1d}|\texttt{02}], [\texttt{1d}|\texttt{11}], [\texttt{1d}|\texttt{13}], [\texttt{1d}|\texttt{21}], [\texttt{1d}|\texttt{23}], [\texttt{1d}|\texttt{31}], [\texttt{1d}|\texttt{33}];$
$\texttt{d} = \texttt{0}, \texttt{1}, \texttt{2}, \texttt{3}\}$

We do not use every cross tile: we will use 155 remaining tiles, $\mathcal{U}_2$, shown at bottom of Figure 18:

$\mathcal{U}_2 :=$ $\{\ [xy|zw]$ with $x = \texttt{3}, \texttt{-3}; y, z, w = \texttt{0}, \texttt{1}, \texttt{2}, \texttt{3}\}$
$\cup\ \{\ [\texttt{-}|zw]$ with $z, w = \texttt{0}, \texttt{1}, \texttt{2}, \texttt{3}\}$
$\cup\ \{\ [xy|\texttt{+}]$ with $xy = \texttt{+}, \texttt{00}, \texttt{01}, \texttt{10}, \texttt{11}, \texttt{12}, \texttt{13}, \texttt{30}, \texttt{-30}, \texttt{32}, \texttt{-32}\} \cup \{[\texttt{-}|\texttt{+}]\}$

The subset $\mathcal{T}_0$ of nineteen $\mathcal{U}_2$ will be common to all of our sets of cross tiles:

$\mathcal{T}_0 :=$ $\{[\texttt{-}|\texttt{+}], [\texttt{00}|\texttt{+}], [\texttt{01}|\texttt{+}], [\texttt{-}|\texttt{00}], [\texttt{-}|\texttt{01}], [\texttt{-}|\texttt{02}], [\texttt{-}|\texttt{03}], [\texttt{-}|\texttt{10}], [\texttt{-}|\texttt{11}], [\texttt{-}|\texttt{12}], [\texttt{-}|\texttt{13}],$
$[\texttt{-}|\texttt{20}], [\texttt{-}|\texttt{21}], [\texttt{-}|\texttt{22}], [\texttt{-}|\texttt{23}], [\texttt{-}|\texttt{30}], [\texttt{-}|\texttt{31}], [\texttt{-}|\texttt{32}], [\texttt{-}|\texttt{33}]\}$

## 4.4 Atomic sets of tiles in $\mathcal{T}_2$

We define atomic sets of tiles, corresponding to the atomic substitutions, and to ordered pairs of atomic substitutions:

For each atomic $\alpha$, we will claim that the sets of tiles in the following table are necessary in order to have a supertile of that type, proven in Section 4.6, and verifiable through the figure below, continuing the construction of Figure 15. (Since it happens that $[\texttt{00}|\texttt{+}]$ and $[\texttt{01}|\texttt{+}]$ must be included in *every* subset of $\mathcal{T}_2$ that tiles the plane, we do not include them in these definitions, but include them instead in the common set of tiles $\mathcal{T}_0$.)

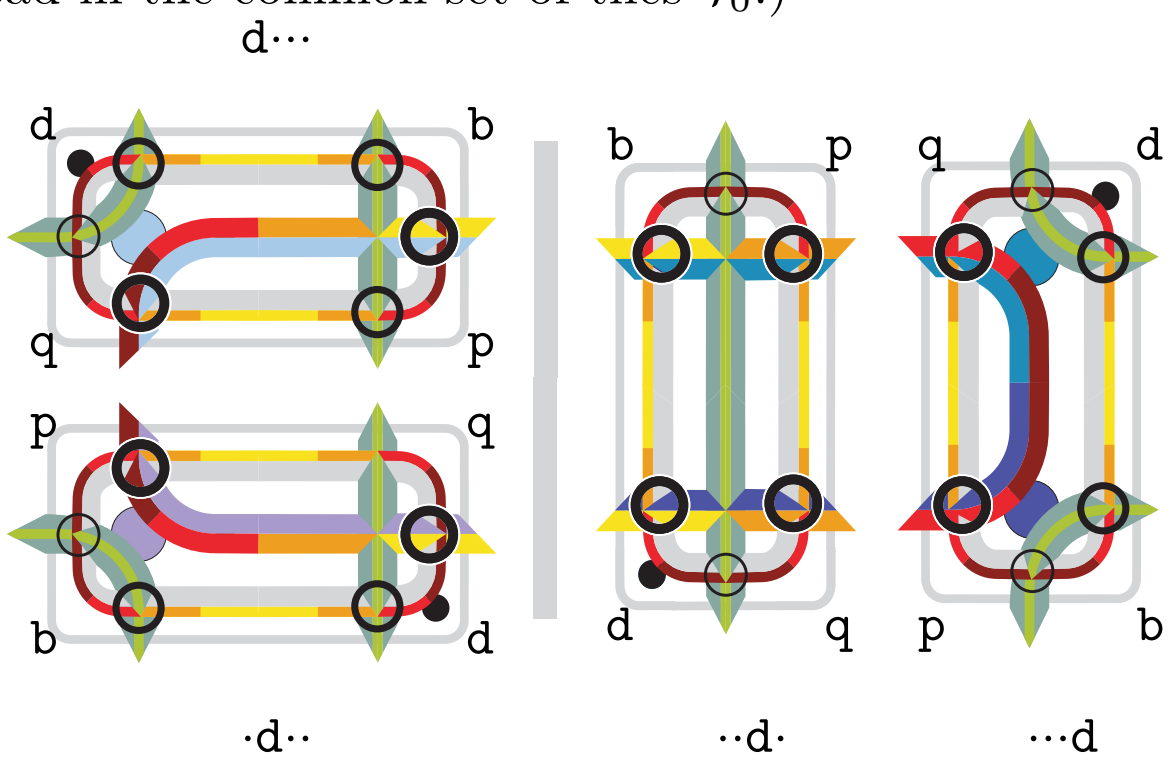

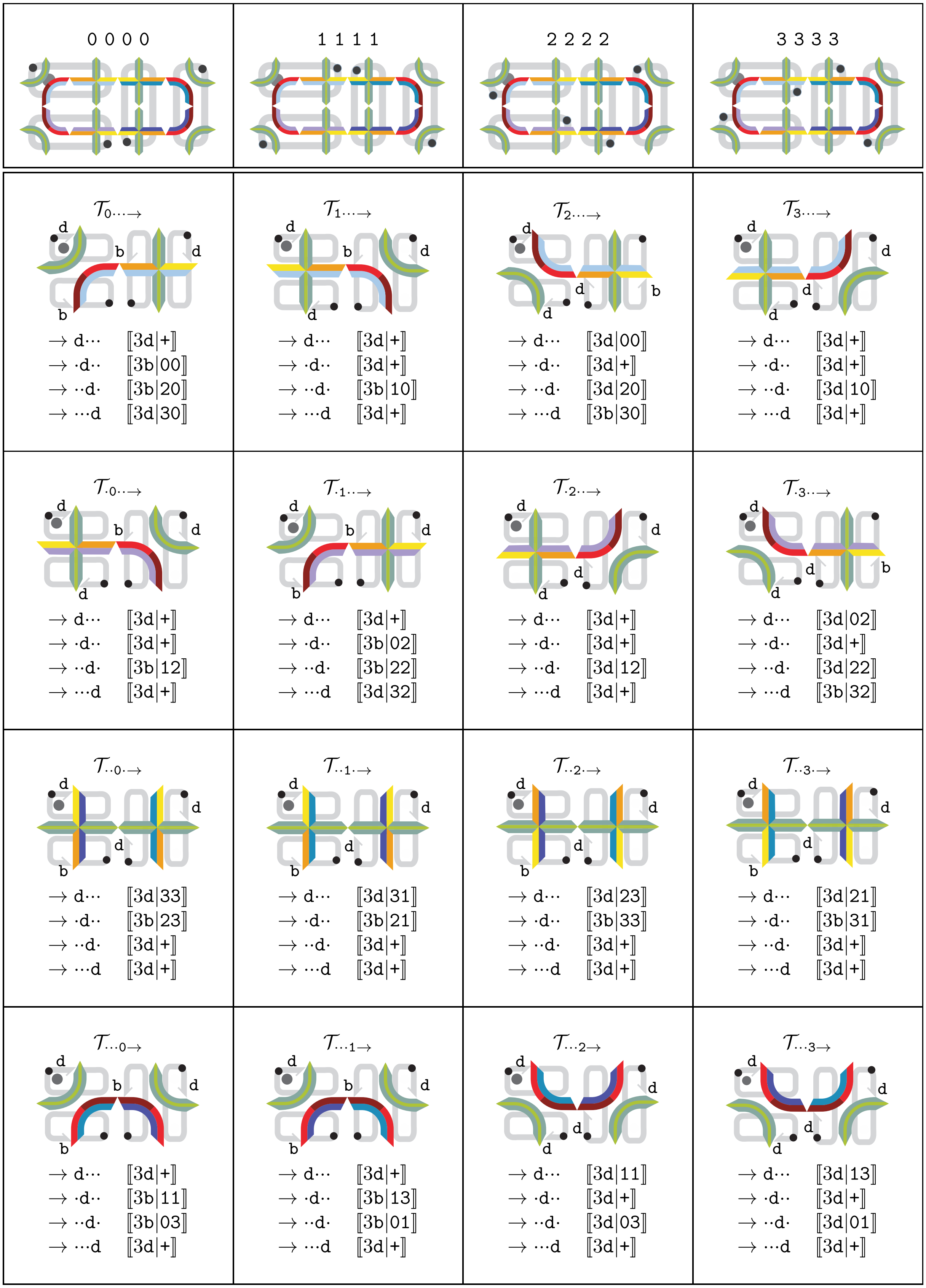

Figure 19: Defining sets of tiles $\mathcal{T}_{\beta\to\alpha}$

| Sets $\mathcal{T}_\alpha$ : | | key tiles | the rest |
|---|---|---|---|
| | $\mathcal{T}_{\mathtt{d}\cdots} =$ | {[1q\|00], [0b\|30], | [1d\|+], [1b\|+], [1p\|+]} |
| | $\mathcal{T}_{\cdot\mathtt{d}\cdot\cdot} =$ | {[0d\|32], [1p\|02], | [1d\|+], [1b\|+], [1q\|+]} |
| | $\mathcal{T}_{\cdot\cdot\mathtt{d}\cdot} =$ | {[1d\|33], [1b\|31], [1q\|23], [1p\|21]} | |
| | $\mathcal{T}_{\cdots\mathtt{d}} =$ | {[1q\|11], [1p\|13], | [1d\|+], [1b\|+]} |

Similarly, in Figure 19, we enumerate, for each ordered pair $\alpha, \beta$ of atomic substitutions, pairs of tiles $\mathcal{T}_{\beta\to\alpha}$, which we claim are necessary for a supertile of type $\alpha$ to have a parent of type $\beta$, as shown in the proof of Theorem 2 by following substitutions on marked $\overline{\overline{\Pi}}$ blocks, and verifiable using the diagrams in the table.

## 4.5 Examples

For example, 1101 has atomic substitutions 1···, ·1··, ··0· and ···1; similarly 1023 has atomic substitutions 1···, ·0··, ··2· and ···3. In addition to $\mathcal{T}_{+_2} \cup \mathcal{T}_0$, $\mathcal{T}_{1101}$ and $\mathcal{T}_{1023}$ include exactly the following atomic sets of tiles:

$\mathcal{T}_{1101} = \mathcal{T}_{+_2} \cup \mathcal{T}_0 \cup$

| | | | |
|---|---|---|---|
| $\mathcal{T}_{1\cdots}$ | [13\|00], [00\|30], [10\|+], [11\|+], [12\|+] | $\mathcal{T}_{\cdot1\cdot\cdot\to\cdot\cdot0\cdot}$ | ⟦31\|21⟧ |
| $\mathcal{T}_{\cdot1\cdot\cdot}$ | [11\|32], [12\|02], [10\|+], [11\|+], [13\|+] | $\mathcal{T}_{\cdot1\cdot\cdot\to\cdots1}$ | ⟦31\|31⟧ |
| $\mathcal{T}_{\cdot\cdot0\cdot}$ | [10\|33], [11\|31], [12\|23], [13\|21] | $\mathcal{T}_{\cdot\cdot0\cdot\to1\cdots}$ | ⟦31\|33⟧ |
| $\mathcal{T}_{\cdots1}$ | [13\|11], [12\|13], [10\|+], [11\|+] | $\mathcal{T}_{\cdot\cdot0\cdot\to\cdot1\cdot\cdot}$ | ⟦30\|23⟧ |
| $\mathcal{T}_{1\cdots\to1\cdots}$ | ⟦31\|+⟧ = ⟦30\|+⟧ | $\mathcal{T}_{\cdot\cdot0\cdot\to\cdot\cdot0\cdot}$ | ⟦30\|+⟧ |
| $\mathcal{T}_{1\cdots\to\cdot1\cdot\cdot}$ | ⟦31\|+⟧ = ⟦30\|+⟧ | $\mathcal{T}_{\cdot\cdot0\cdot\to\cdots1}$ | ⟦30\|+⟧ |
| $\mathcal{T}_{1\cdots\to\cdot\cdot0\cdot}$ | ⟦31\|10⟧ | $\mathcal{T}_{\cdots1\to1\cdots}$ | ⟦30\|+⟧ |
| $\mathcal{T}_{1\cdots\to\cdots1}$ | ⟦30\|+⟧ | $\mathcal{T}_{\cdots1\to\cdot1\cdot\cdot}$ | ⟦30\|13⟧ |
| $\mathcal{T}_{\cdot1\cdot\cdot\to1\cdots}$ | ⟦30\|+⟧ | $\mathcal{T}_{\cdots1\to\cdot\cdot0\cdot}$ | ⟦31\|01⟧ |
| $\mathcal{T}_{\cdot1\cdot\cdot\to\cdot1\cdot\cdot}$ | ⟦30\|01⟧ | $\mathcal{T}_{\cdots1\to\cdots1}$ | ⟦30\|+⟧ |

(32 tiles plus the 16+19 in $\mathcal{T}_{+_2} \cup \mathcal{T}_0$ common to all our aperiodic subsets of $\mathcal{T}_2$.)

$\mathcal{T}_{1023} = \mathcal{T}_{+_2} \cup \mathcal{T}_0 \cup$

| | | | |
|---|---|---|---|
| $\mathcal{T}_{1\cdots}$ | [13\|00], [00\|30], [10\|+], [11\|+], [12\|+] | $\mathcal{T}_{\cdot0\cdot\cdot\to\cdot\cdot2\cdot}$ | ⟦33\|12⟧ |
| $\mathcal{T}_{\cdot0\cdot\cdot}$ | [00\|32], [13\|02], [10\|+], [11\|+], [12\|+] | $\mathcal{T}_{\cdot0\cdot\cdot\to\cdots3}$ | ⟦32\|+⟧ |
| $\mathcal{T}_{\cdot\cdot2\cdot}$ | [12\|33], [13\|31], [10\|23], [11\|21] | $\mathcal{T}_{\cdot\cdot2\cdot\to1\cdots}$ | ⟦31\|23⟧ |
| $\mathcal{T}_{\cdots3}$ | [11\|11], [10\|13], [13\|+], [12\|+] | $\mathcal{T}_{\cdot\cdot2\cdot\to\cdot0\cdot\cdot}$ | ⟦31\|33⟧ |
| $\mathcal{T}_{1\cdots\to1\cdots}$ | ⟦31\|+⟧ = ⟦30\|+⟧ | $\mathcal{T}_{\cdot\cdot2\cdot\to\cdot\cdot2\cdot}$ | ⟦32\|+⟧ |
| $\mathcal{T}_{1\cdots\to\cdot0\cdot\cdot}$ | ⟦30\|+⟧ | $\mathcal{T}_{\cdot\cdot2\cdot\to\cdots3}$ | ⟦32\|+⟧ |
| $\mathcal{T}_{1\cdots\to\cdot\cdot2\cdot}$ | ⟦33\|10⟧ | $\mathcal{T}_{\cdots3\to1\cdots}$ | ⟦31\|13⟧ |
| $\mathcal{T}_{1\cdots\to\cdots3}$ | ⟦33\|+⟧ = ⟦32\|+⟧ | $\mathcal{T}_{\cdots3\to\cdot0\cdot\cdot}$ | ⟦30\|+⟧ |
| $\mathcal{T}_{\cdot0\cdot\cdot\to1\cdots}$ | ⟦30\|+⟧ | $\mathcal{T}_{\cdots3\to\cdot\cdot2\cdot}$ | ⟦32\|01⟧ |
| $\mathcal{T}_{\cdot0\cdot\cdot\to\cdot0\cdot\cdot}$ | ⟦30\|+⟧ | $\mathcal{T}_{\cdots3\to\cdots3}$ | ⟦32\|+⟧ |

(30 tiles plus the 16+19 in $\mathcal{T}_{+_2} \cup \mathcal{T}_0$).

(Sets enforcing non-deterministic substitutions are of course the union of the sets enforcing their constituent deterministic substitution systems, and consequently take longer to write out.)

## 4.6 Proof of Theorem 2

Our main theorem is now well-defined and we turn to its proof. We must show, for each $\overline{\Pi}$ tiling substitution, with symbol $\mathtt{S}$, the set of tiles

$$\mathcal{T}_{\mathtt{S}} := \mathcal{T}_{+_2} \cup \mathcal{T}_0 \cup \left( \bigcup_{\alpha \in \mathtt{S}} \mathcal{T}_\alpha \right) \cup \left( \bigcup_{\alpha,\beta \in \mathtt{S}} \mathcal{T}_{\beta \to \alpha} \right)$$

*does* admit a tiling, and *only* admits tilings that are locally decomposable to $\mathtt{S}$ hierarchical tilings; and moreover that no subset will suffice. The remaining claims will follow from the details of the proof. We have in fact done most of the work, which we mainly need to check in an efficient manner.

We begin by noting there is a natural "forgetting" map from the tiles of $\mathcal{T}_2$ to those of $\mathcal{T}_{\overline{\Pi}} \subset \mathcal{T}_1$: simply ignore the $\mathtt{d}$-channel markings. For the cross tiles, this map takes $[xy|zw] \to [xz]$.

We easily verify that the key tiles $\mathcal{K}_2$ are mapped onto the key tiles in $\mathcal{T}_{\overline{\Pi}}$ (that is, the key tiles $\mathcal{K}_1$ of $\mathcal{T}_1$, with the exception of $[\mathtt{02}]$) and the rest of the tiles $\mathcal{U}_2 \subset \mathcal{T}_2$ to the rest of the tiles $\mathcal{U}_1 \subset \mathcal{T}_1$.

This map gives a local decomposition from tilings and configurations, $\Sigma(\mathcal{T}_2) \subset \mathcal{C}(\mathcal{T}_2)$, admitted by $\mathcal{T}_2$ to $\Sigma(\mathcal{T}_1) \subset \mathcal{C}(\mathcal{T}_1)$. Consequently, any tiling admitted by $\mathcal{T}_2$ (presuming there are any) can only be a (further) marked hierarchical tiling in the {¦|, ≣, $\overline{\Pi}$} tiling substitution.

More particularly, each outward tile must lie in a unique hierarchy of larger and larger $\mathtt{J}$, $\mathtt{H}$ or $\mathtt{U}$ supertiles, themselves framed by outward tiles — cornered outward tiles at the first level only, uncornered outward tiles thereafter.

The definition of our key tiles ensures that these can only be arranged into the $\overline{\Pi}$ tiling substitution, as in Figure 15: Suppose a $\mathtt{J}$ is in the corner $\mathtt{d}$ in a supertile. Then in the supertile we must have tiles of the form $[\mathtt{1}d|\mathtt{0d}]$ and $[\mathtt{0}d|\mathtt{3d}]$ (where $d$ is the orientation of the child, and not important at this moment). In our key tiles, we can have, and only can have the possibility that $\mathtt{d} = \mathtt{0}, \mathtt{2}$.

Similarly, suppose a $\mathtt{U}$ is on the $\mathtt{d}, \mathtt{q}$ end of the supertile. Then in the supertile we must have tiles of the form $[\mathtt{1}d|\mathtt{1d}], [\mathtt{1}q|\mathtt{1q}]$ (where $d, q$ indicate the orientation of the child and are not important at this moment). In our key tiles we can have, and only can have the possibility that $\mathtt{d}, \mathtt{q} = \mathtt{1}, \mathtt{3}$. And for a $\mathtt{H}$ on the $\mathtt{d}, \mathtt{q}$ end of the supertile, we must have $[\mathtt{1}d|\mathtt{2d}], [\mathtt{1}b|\mathtt{2q}], [\mathtt{1}q|\mathtt{3d}], [\mathtt{1}p|\mathtt{3d}]$ — again, we can have and only can have that $\mathtt{d}, \mathtt{q} = \mathtt{1}, \mathtt{3}$

Consequently, any tiling admitted by $\mathcal{T}_2$, if there are any, is a marked $\overline{\Pi}$ hierarchical tiling. We'll distinguish the two $\mathtt{J}$s as $\mathtt{J}_\mathtt{0}$ (in the $\mathtt{0}$ corner of any $\overline{\Pi}$ supertile) and $\mathtt{J}_2$ (in the $\mathtt{2}$ corner). As before, we calculate the blocks and tiles needed to assemble a supertile out of smaller pieces. The analysis is essentially the same:

Let $\mathtt{s}, \mathtt{t}, \mathtt{u}$ and $\mathtt{v}$ be the orientations of $\mathtt{J}_\mathtt{0}, \mathtt{J}_2, \mathtt{H}, \mathtt{U}$ respectively. Recall that each of these acts on the labels *i*, *ii*, ... *ix*, and that a block $\mathtt{X}(x)$ is further specified by the marking $x$ at *ii*.

On a generic supertile $\mathtt{X}(x)$, the children oriented by the domino symmetries $\mathtt{s}, \mathtt{t}, \mathtt{u}, \mathtt{v}$ in $\mathtt{0}, \mathtt{1}, \mathtt{2}, \mathtt{3}$, we have

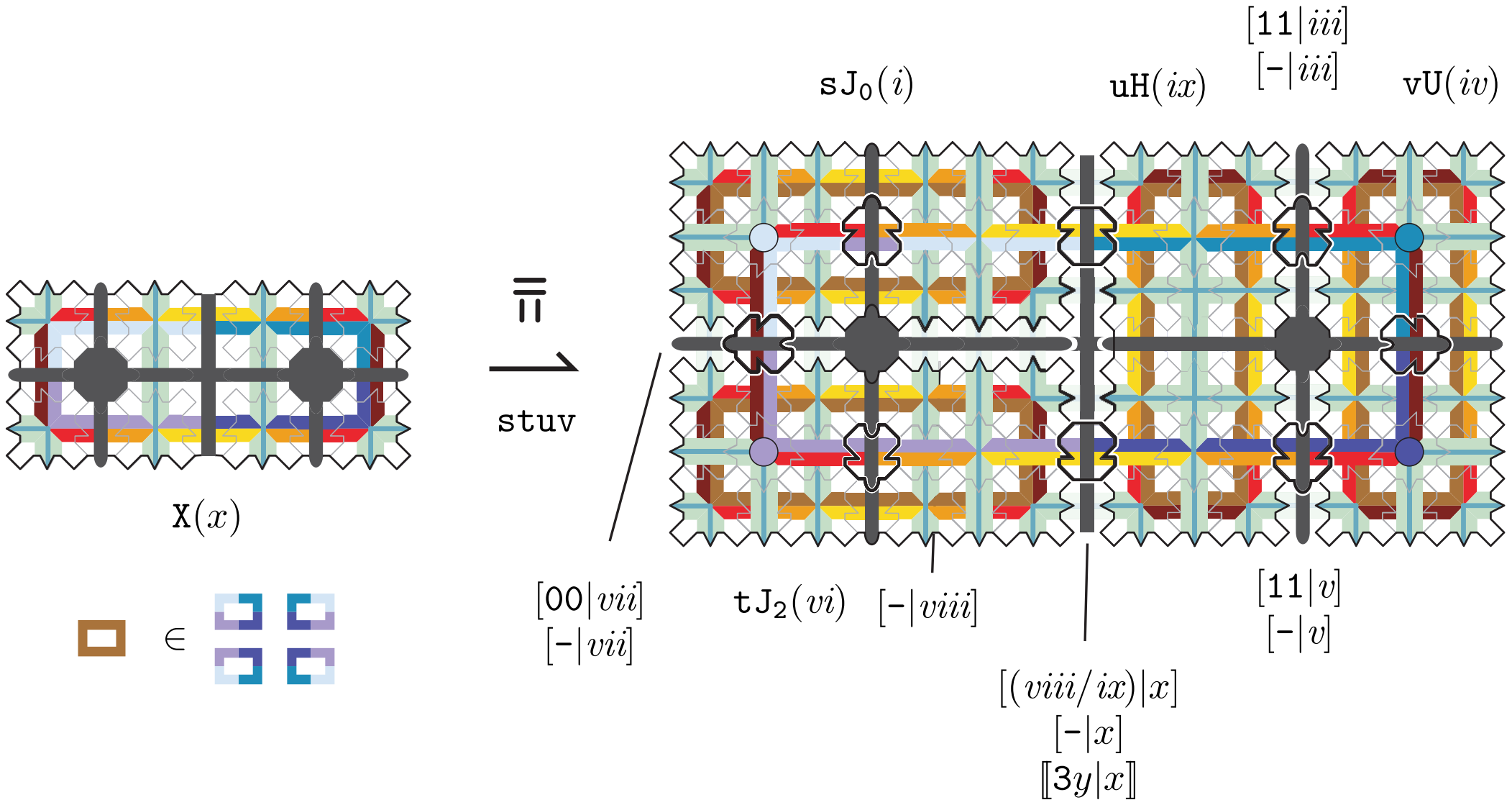

where as before, we resolve $[(viii/ix)|x]$ by examination. We obtain the marking $y$ in the pair $[\![3y|x]\!]$ from the orientation and direction of the marking $x$, that is,

| if $(\mathtt{ab})_x =$ | then $y =$ |
|---|---|
| `+-` or `+0` | 0 |
| `++` | 1 |
| `-+` or `-0` | 2 |
| `--` | 3 |

Substituting a second time on $\mathtt{X}(x)$,

| in the substituted block | we must have the tile |
|---|---|
| $\mathtt{sJ}_0(i)$ | $[\mathtt{10}|i]$ |
| $\mathtt{tJ}_2(vi)$ | $[\mathtt{12}|vi]$ |
| $\mathtt{vU}(iv)$ | $[\mathtt{01}|iv]$ |

In particular, any second level or larger marked $\mathtt{X}(x)$ must include the six tiles $[\mathtt{10}|i]$, $[\mathtt{11}|iii]$, $[\mathtt{01}|iv]$, $[\mathtt{13}|v]$, $[\mathtt{12}|vi]$ and $[\mathtt{00}|vii]$ highlighted in black in the figure above.

Recall, for each atomic symbol $\alpha$ we have defined a set of tiles $\mathcal{T}_\alpha$, and note that these sets are precisely of this form: tiles with horizontal markings $[\mathtt{00}|$, $[\mathtt{01}|$, $[\mathtt{10}|$, $[\mathtt{11}|$, $[\mathtt{12}|$ and $[\mathtt{13}|$ — though we relegate $[\mathtt{00}|+]$, $[\mathtt{01}|+]$ to $\mathcal{T}_0$, as they will be needed in every tiling admitted by any subset of $\mathcal{T}_2$.

We easily check that if $\alpha \in \mathtt{S}$, then any supertile of type $\alpha$ must include the tiles in $\mathcal{T}_\alpha$ — formally we follow the substitutions through and for intuition we examine the figure accompanying the definition of $\mathcal{T}_\alpha$.

Similarly, for each pair of atomic symbols $\alpha, \beta$, we have defined sets of tiles $\mathcal{T}_{\beta\to\alpha}$ and $\mathcal{T}_{\alpha\to\beta}$. Note that each of these is of the form $[\![\mathtt{3}x|yz]\!]$. Following the substitutions above, (examining the tiles outlined in the figure) we see that any supertile of type $\alpha$ within parent supertile of type $\beta$, must include the tiles in $\mathcal{T}_{\beta\to\alpha}$, for intuition examining the figures accompanying the definition above.

In other words, if $\alpha, \beta \in \mathtt{S}$, any set of tiles enforcing the $\mathtt{S}$ tiling substitution must include the sets of tiles $\mathcal{T}_\alpha$, $\mathcal{T}_\beta$, $\mathcal{T}_{\alpha\to\beta}$, $\mathcal{T}_{\beta\to\alpha}$.

Conversely, if a set of tiles $\mathcal{T}$ includes all of

$$\mathcal{T}_{\mathtt{S}} = \mathcal{T}_{+_2} \cup \mathcal{T}_0 \cup \left( \bigcup_{\alpha \in \mathtt{S}} \mathcal{T}_\alpha \right) \cup \left( \bigcup_{\alpha, \beta \in \mathtt{S}} \mathcal{T}_{\beta \to \alpha} \right)$$

then $\mathcal{T}$ *can* form arbitrarily large supertiles arising from the $\mathtt{S}$ tiling substitution system, and hence can tile the plane.

Moreover, noting that every pair of distinct $\mathcal{T}_\alpha$, $\mathcal{T}_\beta$ is mutually exclusive, $\mathcal{T}_{\mathtt{S}}$ can *only* admit supertiles of this form.

Finally, any subset of $\mathcal{T}_2$ that enforces a given $\sigma_{\mathtt{S}}$ must include $\mathcal{T}_{\mathtt{S}}$.

In short, $\mathcal{T}_{\mathtt{S}}$ enforces the $\mathtt{S}$ tiling substitution system, and every subset of $\mathcal{T}_2$ that does as well must contain $\mathcal{T}_{\mathtt{S}}$, concluding the proof of Theorem 2.

# 5 Counting tiling substitution systems and aperiodic sets of tiles

In this section we carefully count out 25,380 hierarchically distinct, regular, elementary ╦╦ tiling substitution systems, sorted into the "lots" of the title of this paper, our tiles and substitutions organized by atomic substitution symbols.

We have increased the number of explicitly described aperiodic sets of tiles hundreds-fold, within a single construction, itself only an example of the kinds of hierarchical control one might routinely, even robustly expect.

This flexibility is easily generalizable, and (given the computational universality at the bottom of this subject) one expects there exist relatively small sets of tiles that have staggeringly complex collections of distinct aperiodic subsets.

But we must take symmetry into account — How many distinct aperiodic sets of tiles have we actually found?

In order to approach this more carefully, we must define what we mean by "distinct" sets of tiles, and by "distinct" tiling substitution systems. Further examples, such as those in Section 6 and [6], may inspire more refined definitions. And the distinction of these, or not, as dynamical systems upon topological spaces awaits industrial scale computation of their invariants, which in most of our cases yet need to be considered.

We will define two sets of tiles, say $\mathcal{T}$ and $\mathcal{T}'$, to be "equivalent" if the sets of tilings $\Sigma(\mathcal{T})$ and $\Sigma(\mathcal{T}')$ are mutually locally decomposible, and "distinct" otherwise.

One might reasonably point out that among our sets of tiles, only 128 (up to mutually local decomposibility) — those enforcing the deterministic ╦╦ tiling substitutions — give tiling spaces that are minimal as dynamical systems. However, the main point of our construction is that there is tremendous flexibility, easily applied and extended if not needlessly restricted. (Each minimal space has only measure zero and lies within the boundary of any of the spaces that properly contain it.)

From any point of view it seems fair to claim that we have substantially increased the number of explicitly described aperiodic sets of tiles.

We will define two tiling substitution systems (here on just the unmarked domino) to be **hierarchically equivalent** if every hierarchical tiling in one system is also a hierarchical tiling in the other, and every supertile (at every level) in the first tiling is also a supertile in the second. More precisely, recall a substitution map $\sigma$ acts on the sets of configurations $\mathcal{C}(\mathcal{T})$ admitted by $\mathcal{T}$; if the map is one-to-one, particularly on tilings, then $\sigma$ has **unique decomposition** and hence has a well-defined inverse on its image. Two tiling substitution systems on the same tiles (here, just the domino) are hierarchically equivalent if and only if this inverse map is equivalent, with the same domain. (The examples below, of non-periodic substitution tilings with non-unique decomposition, are essentially substitution tiling systems that produce the same hierarchical tilings, but are not hierarchically equivalent).

We define an equivalence on our codes that, as we shortly prove, captures hierarchical equivalence on our $\overline{\pi}$ tiling substitions $\sigma_{\mathtt{S}}$ and $\sigma_{\mathtt{S}'}$:

Taking nim addition on the d-channel markings, as in the definition of our code, for $\overline{\pi}$ tiling symbols $\mathtt{S} = \mathtt{stuv}$, $\mathtt{S}' = \mathtt{s}'\mathtt{t}'\mathtt{u}'\mathtt{v}'$, we define $\mathtt{S}' \approx \mathtt{S}$ to hold if and only if $\mathtt{S} = \mathtt{S}'$ or $(\mathtt{stuv})' = (\mathtt{tsuv}) + 3$, that is, if and only if $\mathtt{s} = \mathtt{t}' + 3$, $\mathtt{t} = \mathtt{s}' + 3$, $\mathtt{u} = \mathtt{u}' + 3$ and $\mathtt{v} = \mathtt{v}' + 3$. For example, 1 (02) 3 (012) $\approx$ (13) 2 0 (123) and ·(02)·3 $\approx$ (13)··0.

In a full symbol, each of $\mathtt{s}, \mathtt{t}, \mathtt{u}$ and $\mathtt{v}$ are one of the 15 non-empty subsets of $\{0, 1, 2, 3\}$. There are thus $15^4 = 50625$ distinct full $\overline{\pi}$ tiling substitution symbols, and for each of these we have produced an aperiodic set of tiles. But we have overcounted with respect to this symmetry:

A symbol $\mathtt{S} = \mathtt{stuv}$ is equivalent to one other distinct symbol, unless $\mathtt{s} = \mathtt{t} + 3$, $\mathtt{u} = \mathtt{u} + 3$ and $\mathtt{v} = \mathtt{v} + 3$. There remain fifteen possibilities for $\mathtt{s}$, determining $\mathtt{t}$, and three each for $\mathtt{u}$ and $\mathtt{v}$, namely $\{0, 3\}$, $\{1, 2\}$ and $\{0, 1, 2, 3\}$, or $15 \cdot 3^2 = 135$ symbols not related to any other. All together there are thus $\frac{1}{2}(50625 - 135) + 135 = 25380$ distinct symbols up to this equivalence.

Recall we may add values in $0, 1, 2, 3$ to the d-channel markings, and by extension on tiles, sets of tiles, configurations and tilings: Take $(\mathtt{abcd}) + 2 := (\mathtt{abc}(\mathtt{d} + 2))$ for markings other than (+), leaving (+) alone.

**Proposition 1** *Let* S*,* S′ *be full* $\overline{\pi}$ *tiling symbols.*

1. *The tiling substitutions systems* $\sigma_{\mathtt{S}}$ *and* $\sigma_{\mathtt{S}'}$ *are hierarchically equivalent if and only if* $\mathtt{S} \approx \mathtt{S}'$*;*
2. $\Sigma(\mathcal{T}_{\mathtt{S}})$ *and* $\Sigma(\mathcal{T}_{\mathtt{S}'})$ *are mutually locally decomposible if and only if* $\mathtt{S} \approx \mathtt{S}'$*;*
3. *Moreover, if* $\mathtt{S} \approx \mathtt{S}'$*, then* $\mathcal{T}_{\mathtt{S}} = \mathcal{T}_{\mathtt{S}'} + 2$ *and in fact for each atomic* $\alpha \approx \alpha'$*,* $\beta \approx \beta'$*,* $\mathcal{T}_\alpha = \mathcal{T}'_\alpha + 2$*, and* $\mathcal{T}_{\beta \to \alpha} = \mathcal{T}_{\beta' \to \alpha'} + 2$*.*

**Proof** We first let $\mathtt{S} = \mathtt{stuv}$ and define $\mathtt{S}' = (\mathtt{t} + 3)\ (\mathtt{s} + 3)\ (\mathtt{u} + 3)\ (\mathtt{v} + 3)$; that is, we suppose $\mathtt{S} \approx \mathtt{S}'$. The key observation is that these codes describe precisely the same substitutions, merely with a different framing (switching top and bottom, and the corresponding codes) — this ambiguity in the coding arises from the two-fold-symmetry of the $\overline{\pi}$ configuration.

Illustrating our codes as

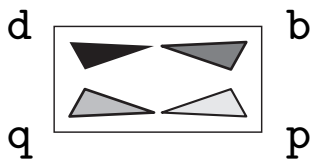

and

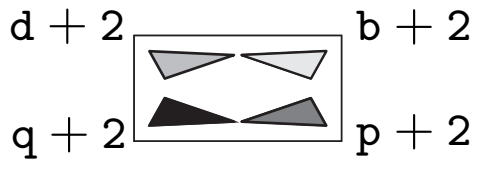

we reflect the substitution rule top-to-bottom, and then replace the codes $0 \leftrightarrow 2$, $1 \leftrightarrow 3$ (that is, applying $+2$), as in the figure below. It is worth checking we can reverse the order: we may just as well interchange the colors and then flip, for the same result:

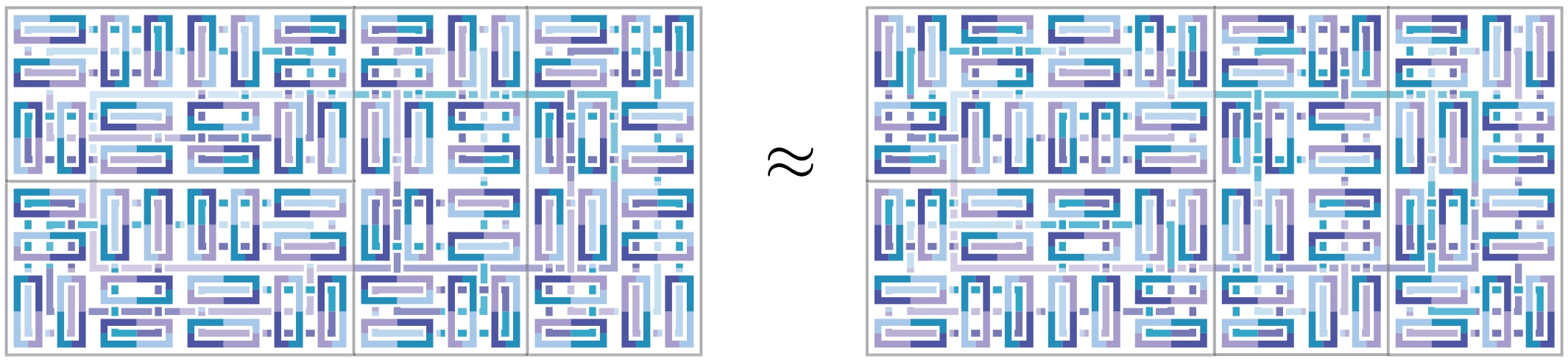

0231 1302

Figure 20: Hierarchies for 0231 $\approx$ 1302. These are congruent if we flip either one vertically and interchange its colors 0 $\leftrightarrow$ 2, 1 $\leftrightarrow$ 3.

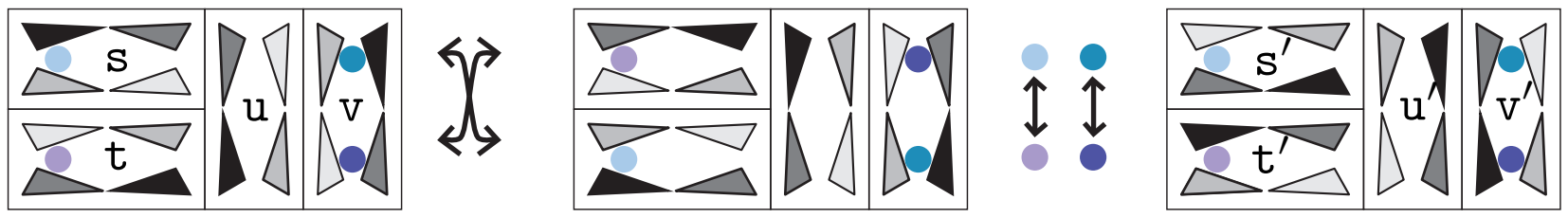

The actual isometries used to define the substitution are the same up to conjugation by this global reflection, and the codes are only for our notational convenience. Up to a global change of coordinates, then, the two substitutions are exactly the same. If $\mathtt{S} \approx \mathtt{S}'$ then $\sigma_{\mathtt{S}}$ and $\sigma_{\mathtt{S}'}$ are hierarchically equivalent.

(Why this specific *ad hoc* relation on symbols? Different encodings will have different two-fold ambiguities, and that of Section 4.2 was chosen as it allows us to define $\mathtt{S} \approx \mathtt{S}'$ fairly simply.)

Because $\mathcal{T}_{\mathtt{S}}$ and $\mathcal{T}_{\mathtt{S}'}$ enforce hierarchically equivalent substitution tiling systems, simply interchanging $\mathtt{0} \leftrightarrow \mathtt{2}$ and $\mathtt{1} \leftrightarrow \mathtt{3}$ in the d-channel (that is, taking each marking $\mathtt{m}$ to $\mathtt{m}+2$) takes any marked supertile of one to a marked supertile of the other. As this is so in either direction, $\Sigma(\mathcal{T}_{\mathtt{S}})$ and $\Sigma(\mathcal{T}_{\mathtt{S}'})$ are mutually locally decomposible. Moreover, since $\mathcal{T}_{\mathtt{S}}$ and $\mathcal{T}_{\mathtt{S}'}$ are precisely the tiles that appear in the tilings of $\Sigma(\mathcal{T}_{\mathtt{S}})$ and $\Sigma(\mathcal{T}_{\mathtt{S}'})$, we have that $\mathcal{T}_{\mathtt{S}} = \mathcal{T}_{\mathtt{S}'} + 2$.

And indeed, this equivalence extends to our atomic sets, as it should: For all atomic $\alpha \approx \alpha'$ and $\beta \approx \beta'$, we claim that $\mathcal{T}_\alpha = \mathcal{T}_{\alpha'} + 2$ and $\mathcal{T}_{\beta\to\alpha} = \mathcal{T}_{\beta'\to\alpha'} + 2$. This can be efficiently checked by hand, from the definitions in Section 4, noticing the following shortcut:

Take atomic $\gamma(\mathtt{d})$ to be any of $\mathtt{d}\cdots$, $\cdot\mathtt{d}\cdot\cdot$, $\cdot\cdot\mathtt{d}\cdot$, $\cdots\mathtt{d}$, and $\gamma'(\mathtt{d})$ to be the corresponding $\cdot\mathtt{d}\cdot\cdot$, $\mathtt{d}\cdots$, $\cdot\cdot\mathtt{d}\cdot$ or $\cdots\mathtt{d}$. In a few minutes we check that for each tile $[xy|zw] \in \mathcal{T}_{\gamma(\mathtt{d})}$ (or $\mathcal{T}_{\beta\to\gamma(\mathtt{d})}$ respectively), the tile $[x(y+1)|w(z+2)] \in \mathcal{T}_{\gamma'(\mathtt{d})}$ (or $\mathcal{T}_{\beta'\to\gamma'(\mathtt{d})}$), (recalling that $[\![\mathtt{31}|+]\!] = [\![\mathtt{30}|+]\!]$ and $[\![\mathtt{33}|+]\!] = [\![\mathtt{32}|+]\!]$, so $[\![\mathtt{3d}|+]\!] = [\![\mathtt{3b}|+]\!]$ and $[\![\mathtt{3q}|+]\!] = [\![\mathtt{3p}|+]\!]$.

We pause for an example; to complete the check does not take much longer. On the left, $\mathcal{T}_{\cdot\mathtt{d}\cdot\cdot} + 2 = \mathcal{T}_{\mathtt{d}\cdots}$ for $\mathtt{d} \in \mathtt{0}, \mathtt{1}, \mathtt{2}, \mathtt{3}$, on the right $\mathcal{T}_{\mathtt{0}\cdots\to\gamma(\mathtt{d})} + 2 = \mathcal{T}_{\cdot\mathtt{3}\cdot\cdot\to\gamma'(\mathtt{d})}$:

| $\mathcal{T}_{\cdot\mathtt{d}\cdot\cdot}$ | | $\mathcal{T}_{\mathtt{d}\cdots}$ |
|---|---|---|
| [0d\|32] | [0(d + 1)\|3(2 + 2)] | [0b\|30] |
| [1p\|02] | [1(p + 1)\|0(2 + 2)] | [1q\|00] |
| [1d\|+] | [1(d + 1)\|+] | [1b\|+] |
| [1b\|+] | [1(b + 1)\|+] | [1d\|+] |
| [1q\|+] | [1(q + 1)\|+] | [1p\|+] |

| $\mathcal{T}_{\mathtt{0}\cdots\to}$ | | $\mathcal{T}_{\cdot\mathtt{3}\cdot\cdot\to}$ | |
|---|---|---|---|
| d··· | ⟦3d\|+⟧ | ·d·· | ⟦3b\|+⟧ = ⟦3d\|+⟧ |
| ·d·· | ⟦3b\|00⟧ | d··· | ⟦3d\|02⟧ |
| ··d· | ⟦3b\|20⟧ | ··d· | ⟦3d\|22⟧ |
| ···d | ⟦3d\|30⟧ | ···d | ⟦3b\|32⟧ |

Once completed, the full check verifies the claim: $[xy|zw] \in \mathcal{T}_{\gamma'(\mathtt{d})}$ (or $\mathcal{T}_{\beta'\to\gamma'(\mathtt{d})}$) if and only if $[x(y+3)|zw] \in \mathcal{T}_{\gamma'(\mathtt{d}+3)}$ (or $\mathcal{T}_{\beta\to\gamma(\mathtt{d}+3)}$). If $\alpha = \gamma(\mathtt{d})$, then $\alpha' = \gamma'(\mathtt{d}+3)$ and so $[xy|zw] \in \mathcal{T}_\alpha$ (or $\mathcal{T}_{\beta\to\alpha}$) if and only if $[x(y+2)|z(w+2)] \in \mathcal{T}_{\alpha'}$ (or $\mathcal{T}_{\beta'\to\alpha'}$).

Next, conversely, suppose that $\mathtt{S} \not\approx \mathtt{S}'$. We first establish that $\sigma_{\mathtt{S}}$ and $\sigma_{\mathtt{S}'}$ are not hierarchically equivalent. Because $\mathtt{S} \not\approx \mathtt{S}'$, there must be a deterministic symbol $\mathtt{S}_0$ in, say $\mathtt{S}$ such that $\mathtt{S}_0' \notin \mathtt{S}'$. Consequently, there are deterministic $\mathtt{S}_0 = \mathtt{s}_0\mathtt{t}_0\mathtt{u}_0\mathtt{v}_0 \not\approx \mathtt{S}_1 = \mathtt{s}_1\mathtt{t}_1\mathtt{u}_1\mathtt{v}_1$, with $\mathtt{S}_0 \in \mathtt{S}$, $\mathtt{S}_1 \in \mathtt{S}'$.

For some $\mathtt{x} \in \{\mathtt{s}, \mathtt{t}, \mathtt{u}, \mathtt{v}\}$, $\mathtt{x}_0 \neq \mathtt{x}_1$. Suppose first that $\mathtt{x}_0 = \mathtt{x}_1 + \mathtt{1}$ or $\mathtt{x}_1 + \mathtt{3}$ (This case with $\mathtt{x} = \mathtt{s}$ is illustrated at left below). Then every 2-level supertile in $\sigma_{\mathtt{S}_0}$ is distinct from every 2-level supertile in $\sigma_{\mathtt{S}_1}$ — the 1-level supertiles in the $\mathtt{x}$ positions are not in the same orientation. As this is preserved under substitution, for every $n \geq 2$, the corresponding $n$-level supertiles in $\sigma_{\mathtt{S}_0}$, $\sigma_{\mathtt{S}_1}$ are distinct.

Suppose on the other hand that $\mathtt{x}_0 = \mathtt{x}_1 + \mathtt{2}$. (This cas with $\mathtt{x} = \mathtt{s}$ is illustrated at center below, and with $\mathtt{x} = \mathtt{u}$ at right. Each subcase is similar, as we leave to the reader to more fully verify.) Every 3-level supertile in $\sigma_{\mathtt{x}_0}$ is distinct from every 3-level supertile in $\sigma_{\mathtt{x}_1}$ — the 1-level supertiles in the $\mathtt{u}$ positions within the 2-level supertiles in the $\mathtt{x}$ position are are not the same, (keeping in mind that $\mathtt{u}_0 = \mathtt{u}_1$ or $\mathtt{u}_1 + \mathtt{2}$). As this is preserved under substitution, for every $n \geq 3$, the corresponding $n$-level supertiles in $\sigma_{\mathtt{S}_0}$, $\sigma_{\mathtt{S}_1}$ are distinct.

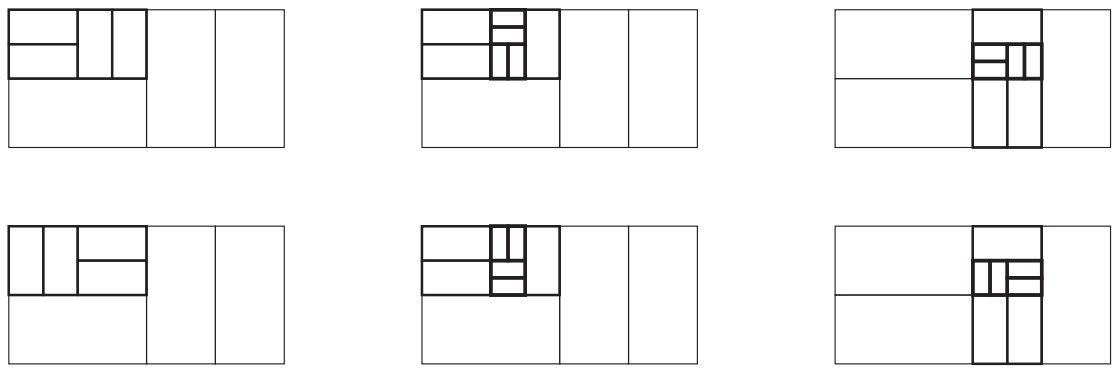

In short, if $\mathtt{S} \not\approx \mathtt{S}'$ then $\sigma_{\mathtt{S}}$ and $\sigma_{\mathtt{S}'}$ are not hierarchically equivalent. Finally we prove that if $\mathtt{S} \not\approx \mathtt{S}'$, then $\Sigma(\mathcal{T}_{\mathtt{S}})$ and $\Sigma(\mathcal{T}_{\mathtt{S}'})$ are not mutually locally decomposable.

It suffices to assume that $\mathtt{S}$ and $\mathtt{S}'$ are deterministic, since for each deterministic $\mathtt{S}_0 \subset \mathtt{S}$ and $\mathtt{S}_0' \subset \mathtt{S}'$, $\Sigma(\mathcal{T}_{\mathtt{S}_0}) \subset \Sigma(\mathcal{T}_{\mathtt{S}})$ and $\Sigma(\mathcal{T}_{\mathtt{S}_0'}) \subset \Sigma(\mathcal{T}_{\mathtt{S}'})$. Any local decomposition from one of the larger sets to the other would restrict to a local decomposition on the smaller sets as well; consequently if for no pair of smaller sets there is a local decomposition (noting they are minimal as subshifts), neither is there among the larger ones. So without loss of generality, we assume $\mathtt{S}$ and $\mathtt{S}'$ are deterministic, and are each specified by a single substitution rule.

Consider any radius $R > 0$, any tiling $T \in \Sigma(\mathcal{T}_{\mathtt{S}})$ and any $T' \in \Sigma(\mathcal{T}_{\mathtt{S}'})$. We will show there exists a pair of distinct tiles $t_1, t_2 \in T'$ so that there exist congruent neighborhoods $C_1$, $C_2$ in $T$, containing all tiles within $R$ of the support of $t_1$, $t_2$ respectively. Consequently, there can be no well-defined local decompostion taking $C_1$, $C_2$ to configurations containing $t_1$ or $t_2$. Hence there can be no well-defined local map on neighborhoods of radius $R$ from $T$ to $T'$, and so no local decomposition from $\Sigma(\mathcal{T}_{\mathtt{S}})$ to $\Sigma(\mathcal{T}_{\mathtt{S}'})$. We will have shown that $\mathcal{T}_{\mathtt{S}}$ and $\mathcal{T}_{\mathtt{S}'}$ are not mutually locally decomposible.

The substitutions have well-defined inverses (deflations), $\sigma_{\mathtt{S}}^{-1}$, $\sigma_{\mathtt{S}'}^{-1}$ on the tilings in $\Sigma(\mathcal{T}_{\mathtt{S}}), \Sigma(\mathcal{T}_{\mathtt{S}'})$. For some $n \in \mathbb{N}$, $2^n > R$ and consider $T_0 = \sigma_{\mathtt{S}}^{-n}(T)$ and $T_0' = \sigma_{\mathtt{S}'}^{-n}(T')$.

Because the substitutions $\sigma_{\mathtt{S}}$ and $\sigma_{\mathtt{S}'}$ are not hierarchically equivalent, there exist some pair of locations $x, g_0 x \in \mathbb{E}^2$, $g_0$ an isometry, such that $x$, $g_0 x$ are in the interiors of congruent marked blocks $C_0$, $g_0 C_0$ in $T_0$ and $x$ and $g_0 x$ are on non-congruent cornered tiles in the interior marked blocks in $T_0'$.

Let $g = 2^n g_0 2^{-n}$. Let $C_1$ be the minimal configuration in $T$ containing $2^n[\![C_0]\!]$; then $C_2 = gC$ contains $2^n[\![g_0 C_0]\!]$ and morever $gC \in T$. But cornered tiles are preserved under substitution. The tiles $t_1, t_2$ in $T'$ at $2^n x$ and $2^n g_0 x = g(2^n x)$ are not congruent, and we have completed the proof.

ddbb

0011 $\approx$ 3322

1100 $\approx$ 2233

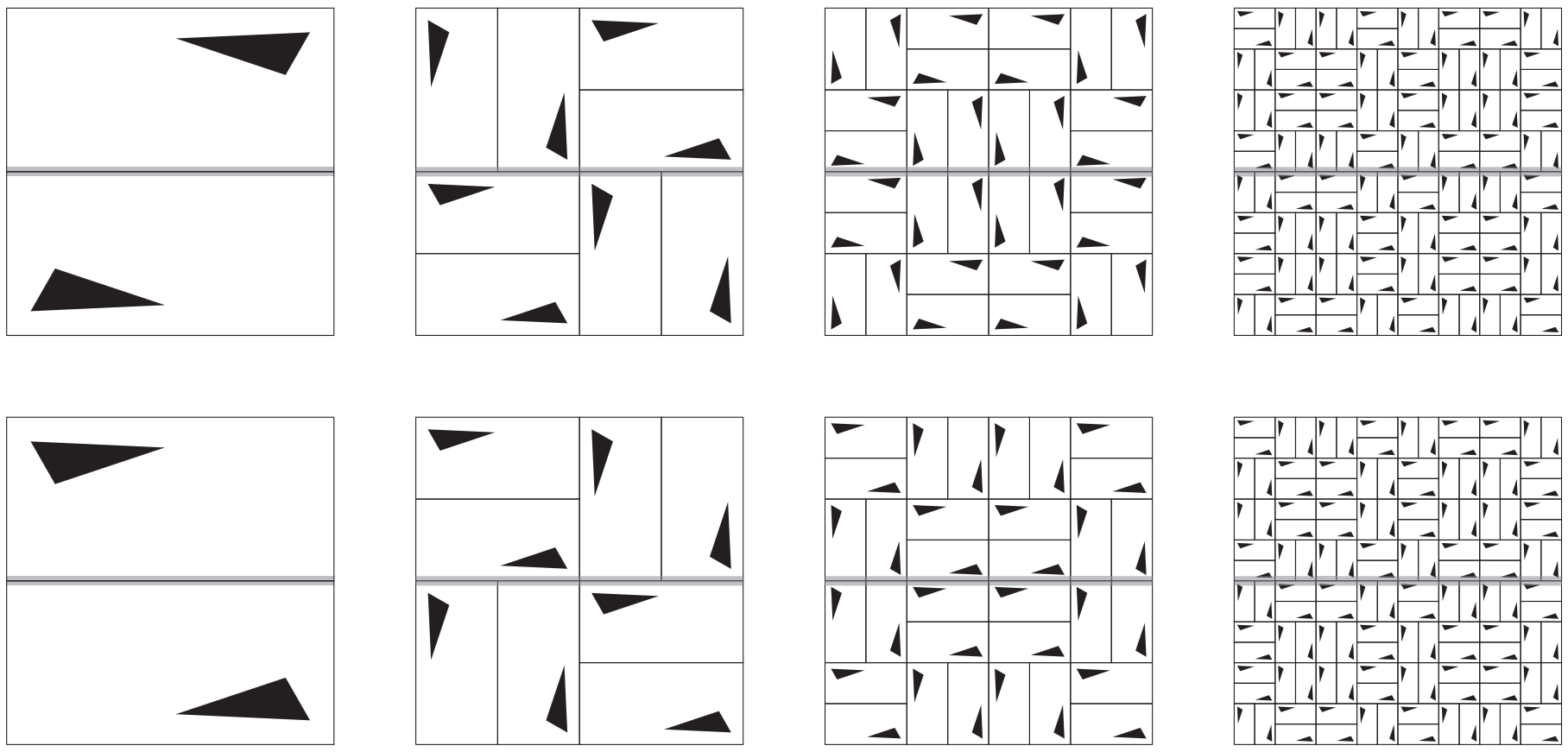

Figure 21: The 0011 $\approx$ 3322 and 1100 $\approx$ 2233 (all of the form ddbb) substitutions are not hierarchically equivalent, but have equivalent hierarchical tilings (note that the two configurations at right in the figure are the same) — the substitution tilings from the two substitutions are the same, and each substitution tiling thus has two distinct hierarchies that can be imposed upon it. The substitutions share the same substitution tiling space, but with distinct actions upon it. Interestingly, as Harriss points out, these two systems factor onto the product of Thue-Morse symbolic substitutions. The substitutions of the form dddd have non-unique hierarchy but are periodic. There are two further deterministic non-periodic, non-unique decomposition examples, ddpp, ddqq, and we may take unions of these. Further comments appear in [6].

# 6 Non-periodic substitution tilings with non-unique decomposition

Every substitution tiling system with unique decomposition contains only non-periodic substitution tilings — from [2] onwards, this observation has been at the heart of every construction of aperiodic sets of tiles that enforce hierarchical tilings.

Moreover, Mossé [11] proved the converse in symbolic substitution systems (essentially equivalent to tiling substitution systems on the line): if the subshift arising from a symbolic substitution system is aperiodic then the system has unique decomposition. Solomyak [18] extended this to two-dimensional subshifts, that is, substitution tilings on which the isometries are all translations. Hunton, Radin and Sadun [9] further showed that on substitution tiling spaces in which the tiles are oriented only finitely, or densely, if non-periodic then the substitution action on the substitution tiling space is one-to-one, i.e. each tiling in the space has a unique hierarchy under that substitution. (They sketch a higher dimensional example of a non-periodic “quasi-substitution” that has non-one-to-one action).

However, the tilings in this section are non-periodic yet have non-unique hierarchy, with a one-to-one substitution action on the tiling space! The trick is they have another, distinct one-to-one action as well, that is we give pairs of distinct substitution systems that have the same underlying substitution tiling space, yet have distinct substitution actions upon it.

(Our tilings fall under Solomyak’s theorem, in any case, if we modify our dominos so that they do not need to be rotated or reflected to carry out the substitutions, obtaining distinct tiling spaces; but this requires multiple types of marked domino. However our substitutions are defined with *unmarked* tiles.)

**Proposition 2** *Let* $\mathtt{S} = \mathtt{stuv}$ *be a full* $\overline{\sqcap\sqcap}$ *symbol, with* $\mathtt{s}, \mathtt{t}, \mathtt{u}, \mathtt{v} \in \{0, 1, 2, 3, (02), (13)\}$ *and* $\mathtt{s} = \mathtt{t}$, $\mathtt{u} = \mathtt{v}$. *(Note that* $\mathtt{S} \approx \mathtt{S} + 2$*.) Let* $\mathtt{S}' = \mathtt{S} + 1$ *or* $\mathtt{S} + 3$. *Although* $\sigma_{\mathtt{S}}$ *and* $\sigma_{\mathtt{S}'}$ *are* *not hierarchically equivalent,* every *hierarchical tiling of one is a hierarchical tiling of the other. Consequently, every hierarchical tiling of either has non-unique decomposition. Moreover if* $\mathtt{s} \neq \mathtt{u}$*, then every hierarchical tiling of either is non-periodic.*

There are thus twelve quadruples of substitution systems with equal hierarchical tilings, each quadruple in two hierarchically equivalent pairs: $(\mathtt{d}^2, (\mathtt{dq})^2) \times ((\mathtt{d}, \mathtt{b\,q}, \mathtt{p}, (\mathtt{dq}), (\mathtt{bp}))$. Four of these quadruples are deterministic.

**Proof** We have already established that any $\mathtt{S} \not\approx \mathtt{S}+\mathtt{1}$, $\mathtt{S} \not\approx \mathtt{S}+\mathtt{3}$. Take $\mathtt{S}$, $\mathtt{S}'$ as in the theorem. Because $\mathtt{s} = \mathtt{t}$, $\mathtt{u} = \mathtt{v}$, by induction each supertile consists of two rotationally symmetric squares. However, for level $n > 1$, these squares can themselves be partitioned into supertiles in two different ways, precisely corresponding to the substitutions $\mathtt{S}$, $\mathtt{S}'$. That the tilings are non-periodic can be verified easily by considering marked dominos — the orientations of the tiles are visible and imply unique decomposition, hence are non-periodic. Our *unmarked* dominos have substitution tilings with *non-unique* decomposition. But because they are just marked dominos with the markings erased as a local derivation, the unmarked domino substitution tilings are non-periodic.

In fact, with this example in front of us, we see this phenomenon already occurs in simpler settings, and possibly quite widely: The one dimensional non-periodic fibonacci symbolic substitution system has non-unique decomposition if we allow ourselves to change the global orientation of the hierarchy. Harriss points out that the Rauzy tiles also have a natural two fold symmetry in their substitutions and non-periodic substitution tilings with non-unique decomposition.